\documentclass[hidelinks,onefignum,onetabnum]{siamart251216}
\newsiamremark{assumption}{Assumption}
\crefname{assumption}{Assumption}{Assumptions}

\usepackage{lipsum}
\usepackage[varbb]{newpxmath}
\usepackage{amsfonts,amsmath}
\usepackage{graphicx}
\usepackage{epstopdf}
\usepackage{algorithmic}
\ifpdf
\DeclareGraphicsExtensions{.eps,.pdf,.png,.jpg}
\else
\DeclareGraphicsExtensions{.eps}
\fi

\newsiamremark{remark}{Remark}
\newsiamremark{hypothesis}{Hypothesis}
\crefname{hypothesis}{Hypothesis}{Hypotheses}
\newsiamthm{claim}{Claim}
\newsiamremark{fact}{Fact}
\crefname{fact}{Fact}{Facts}

\headers{Optimal control of nonmonotone singular equations}{B. Vexler}

\title{Optimal Control of Nonmonotone Semilinear Elliptic Equations with Possibly Singular Linearizations 
}

\author{Boris Vexler\thanks{Technical University of Munich, Germany
		(\email{vexler@tum.de},
		\url{https://www.math.cit.tum.de/math/vexler});
		Kutaisi International University, Georgia.}}
	
\usepackage{amsopn}

\ifpdf
\hypersetup{
  pdftitle={Optimal control},
  pdfauthor={B. Vexler}
}
\fi

\newcommand{\R}{\varmathbb{R}}

\newcommand{\norm}[1]{\lVert#1\rVert}

\newcommand{\ltwonorm}[1]{\norm{#1}_{L^2(\Omega)}}

\newcommand{\half}{\frac{1}{2}}

\newcommand{\oq}{\bar q}

\renewcommand{\phi}{\varphi}

\begin{document}

\maketitle

\begin{abstract}
We study optimal control problems governed by nonmonotone
semilinear elliptic equations with pointwise control constraints.
The state equation need not define a single-valued
control-to-state mapping, and its linearization may be singular.
We prove the existence of globally optimal control-state pairs
whenever the admissible set is nonempty.
Under the assumption that the inactive control set has positive
Lebesgue measure, a penalization argument combined with unique
continuation yields normal first-order optimality conditions
and uniqueness of the adjoint satisfying the full optimality
system.
We further establish second-order necessary and sufficient
conditions of no-gap type.
\end{abstract}



\section{Introduction}
In this paper, we study optimal control problems governed
by nonmonotone semilinear elliptic equations with pointwise
control constraints.
The state equation need not define a single-valued
control-to-state mapping, and its linearization may be singular.
We formulate the problem in terms of admissible control-state
pairs and develop normal first-order optimality conditions
and a second-order optimality theory.

As a model problem, we study the following optimal control problem
\begin{equation}
\text{Minimize } J(q,u) = \half\ltwonorm{u-u_d}^2 +\frac{\alpha}{2} \norm{q}_{L^2(\omega)}^2
\end{equation}
subject to
\begin{equation}
\begin{aligned}
	-\Delta u +d(u) &=&&  q \chi_\omega & \text{in} && \Omega\\
	              u &=&& 0              &\text{on}  && \partial \Omega
\end{aligned}
\end{equation}
and pointwise control constraints
\begin{equation}
q_a \le q(x) \le q_b \quad \text{for almost all } x \in \omega.
\end{equation}
The precise assumptions for this problem setting are discussed below. In particular, we will neither assume monotonicity nor (upper) growth conditions on the nonlinearity $d$. We will essentially assume $d \in C^1(\R)$ and $d'(t) \ge -b$ for all $t \in \R$ with some $b\ge 0$, see examples below. For second order analysis additional assumptions on the second derivative of $d$ are required.

In the majority of published results on optimal control of semilinear elliptic equations, see, e.g., the monographs \cite{Troeltzsch:2010,HinzePinnauUlbrichUlbrich:2009,VexlerMeidner:2025}, the nonlinearity is assumed to be monotone requiring $d'(t)\ge 0$. In this case, the state equation admits a unique solution for every control, while the linearized state operator is coercive and hence continuously invertible. Consequently, the control-to-state operator is well defined and differentiable, and the adjoint equation can be introduced in the usual way.

First- and second-order optimality conditions for a class
of nonmonotone and noncoercive semilinear elliptic operators
were established by Casas, Mateos, and R\"osch
\cite{CasasMateosRoesch:2020}.
In their setting, nonmonotonicity arises from a convection
term.
The state equation nevertheless admits a unique solution
for every control, and the linearized state operator is
invertible. These properties provide a twice continuously differentiable
control-to-state mapping and allow a second-order analysis
of the reduced objective functional.

Optimal control problems for which the state equation does
not define a single-valued control-to-state mapping have
been studied since the early work of Lions \cite{Lions:1985}.
A further important contribution is due
to Bonnans and Casas \cite{BonnansCasas:1989}, who considered
semilinear multistate systems with control and state
constraints. A central difficulty in such problems is the derivation
of necessary optimality conditions.
For general control constraints, these conditions are
often obtained in Fritz--John form, allowing the multiplier
associated with the objective functional to vanish.
Excluding this abnormal case requires additional arguments
and assumptions on the problem. Normal optimality conditions have been established
for several particular settings.

Bonnans and Casas \cite[Theorem~5]{BonnansCasas:1989}
use unique continuation to obtain normality under
additional assumptions on the admissible controls,
including the availability of unrestricted control
variations on an open subset of the domain.
For an ill-posed elliptic equation with an exponential
nonlinearity, Casas, Kavian, and Puel
\cite{CasasKavianPuel1998} derive first-order conditions
by penalization and localization.
Under additional regularity assumptions, they obtain
normality when the admissible control set permits
localized two-sided perturbations on an open subset;
see \cite[Corollary~5.5]{CasasKavianPuel1998}. Fursikov \cite[Chapter~2, Section~7.1]{Fursikov2000}
develops a Lagrange multiplier approach to nonlinear
elliptic control problems that allows a nontrivial kernel
of the linearized state operator.
He obtains normal conditions for unconstrained distributed
controls and identifies further cases with normality,
including problems with a one-sided control constraint
and an additional structural condition on the nonlinearity.
Further first-order results for non-well-posed elliptic
control problems are given by Wang and Wang
\cite{WangWang:2003}.
Related problems governed by $p$-Laplacian equations
are studied by Lou and Luan \cite{LouLuan:2024}. 

Second-order necessary and sufficient optimality conditions
in a setting without a well-defined control-to-state mapping
were obtained by Liu and Wang \cite{LiuWang:2021} for
semilinear parabolic equations with a periodicity constraint.
Their results concern unconstrained controls and use
controllability arguments to construct suitable linearized
state-control directions.

For the class of problems described above, our main
contributions are as follows:

\begin{itemize}
	\item
	We establish the existence of a globally optimal
	control-state pair for the considered class of problems;
	see Theorem~\ref{thm:existence}.
	
	\item
	We derive a normal first-order optimality system
	with the same structure as in problems with
	well-posed state equations;
	see Theorem~\ref{thm:optimality_system}.
	The only additional structural assumption on the
	local solution is that the set on which the control
	constraints are inactive has positive Lebesgue measure;
	see Assumption~\ref{ass:inactive_set}.
	Moreover, the adjoint state satisfying the full
	optimality system is unique, even though the adjoint
	equation alone may admit multiple solutions.
	We also present an example showing that normality
	can fail when the inactive-set assumption is omitted.
	
	\item
    We establish a second-order
	sufficient optimality condition in terms of strict
	positivity of the second variation on all nonzero
	critical control-state directions;
	see Theorem~\ref{thm:second_order_sufficient}.
	For an admissible pair satisfying the first-order
	optimality system, this condition yields quadratic
	growth in the joint
	$L^2(\omega)\times L^2(\Omega)$ topology and hence
	strict local optimality.
	
	\item
	We complement the sufficient condition with a
	second-order necessary optimality condition;
	see Theorem~\ref{thm:second_order_necessary}.
	This result
	requires the same structural Assumption
	\ref{ass:inactive_set} on the local solution as the
	first-order theory. The necessary and sufficient conditions differ only in nonnegativity
	versus strict positivity on nonzero critical directions,
	providing second-order conditions of no-gap type.
\end{itemize}

The proof of the first-order optimality system combines
a penalty approach with an appropriate unique continuation
principle.
The penalized formulation may also serve as a basis for
homotopy-type optimization algorithms, in which a sequence
of optimal control problems with well-posed state equations
is solved approximately. The development and analysis of such algorithms are
the subject of ongoing research.

\section{Problem Setting and Auxiliary Results}
We consider the following optimal control problem:
\begin{subequations}\label{eq:ocp}
	\begin{equation}\label{eq:ocp_cost}
		\text{Minimize}\quad
		J(q,u)
		=
		\frac12\|u-u_d\|_{L^2(\Omega)}^2
		+
		\frac{\alpha}{2}\|q\|_{L^2(\omega)}^2
	\end{equation}
	subject to the state equation
	\begin{equation}\label{eq:ocp_state}
		\begin{aligned}
			-\Delta u+d(x,u) &= q\chi_\omega
			&& \text{in }\Omega,\\
			u &= 0
			&& \text{on }\partial\Omega,
		\end{aligned}
	\end{equation}
	and the control constraints
	\begin{equation}\label{eq:ocp_constraints}
		q\in Q_{\mathrm{ad}}
		:=
		\left\{
		q\in L^2(\omega):
		q_a\le q(x)\le q_b
		\text{ for almost all }x\in\omega
		\right\}.
	\end{equation}
\end{subequations}
Here, $q\chi_\omega$ denotes the extension of $q$ by zero to $\Omega$.

Throughout this paper, we impose the following assumptions.
The set $\Omega\subset\mathbb{R}^N$, $N\in\{2,3\}$, is a bounded,
connected, open set with Lipschitz boundary, and
$\omega\subset\Omega$ is a measurable subset of positive Lebesgue measure.
The desired state satisfies $u_d\in L^2(\Omega)$, and $\alpha>0$.
The control bounds are constants
\[
q_a\in\mathbb{R}\cup\{-\infty\},
\qquad
q_b\in\mathbb{R}\cup\{+\infty\},
\qquad q_a<q_b.
\]
An infinite bound means that the corresponding one-sided constraint
is absent.

We impose the following assumption on the nonlinearity
throughout the paper.

\begin{assumption}\label{ass:d}
	The nonlinearity $d:\Omega\times\mathbb R\to\mathbb R$
	satisfies the following conditions:
	\begin{itemize}
		\item For every $t\in\mathbb R$, the mapping
		$x\mapsto d(x,t)$ is measurable, and
		$d(\cdot,0)\in L^2(\Omega)$.
		
		\item For almost every $x\in\Omega$, the mapping
		$t\mapsto d(x,t)$ is continuously differentiable.
		
		\item For every $M>0$, there exists a constant
		$c_M>0$ such that
		\begin{equation}\label{eq:d_derivative_bound}
			|\partial_t d(x,t)|\le c_M
			\qquad
			\text{for almost every }x\in\Omega
			\text{ and all }|t|\le M.
		\end{equation}
		
		\item There exists a constant $b\ge0$ such that
		\begin{equation}\label{eq:d_lower_bound}
			\partial_t d(x,t)\ge-b
			\qquad
			\text{for almost every }x\in\Omega
			\text{ and all }t\in\mathbb R.
		\end{equation}
	\end{itemize}
\end{assumption}

Assumption~\ref{ass:d} implies
$d(\cdot,t)\in L^2(\Omega)$ for every $t\in\mathbb R$.
We impose neither monotonicity of $d(x,\cdot)$ nor
growth conditions as $|t|\to\infty$.

For the second-order analysis below, we additionally impose
the following regularity assumption.

\begin{assumption}\label{ass:d_second_order}
	For almost every $x\in\Omega$, the mapping
	$t\mapsto d(x,t)$ is twice continuously differentiable.
	Moreover, for every $M>0$, there exists a constant
	$L_M>0$ such that
	\begin{equation}\label{eq:d_second_derivative_lipschitz}
		|\partial_{tt}d(x,t_1)-\partial_{tt}d(x,t_2)|
		\le L_M|t_1-t_2|
	\end{equation}
	for almost every $x\in\Omega$ and all
	$t_1,t_2\in[-M,M]$.
\end{assumption}
Together with Assumption~\ref{ass:d},
Assumption~\ref{ass:d_second_order} implies that
$\partial_{tt}d$ is bounded on bounded $t$-intervals,
uniformly with respect to almost every $x\in\Omega$.

Typical examples of nonlinearities satisfying these assumptions are
\[
d(x,t)=\mu(t^3-t), \qquad
d(x,t)=t^5-2t^3, \qquad
d(x,t)=e^t-\mu t,
\]
where $\mu>0$ may be arbitrarily large.
More generally, every polynomial in $t$ of odd degree with a positive
leading coefficient satisfies the above assumptions.

An example with explicit spatial dependence is
\[
d(x,t)=a(x)(t^3-t)+c(x),
\]
where $a\in L^\infty(\Omega)$ is nonnegative almost everywhere and
$c\in L^2(\Omega)$. Indeed,
\[
\partial_t d(x,t)=a(x)(3t^2-1)
\ge -\|a\|_{L^\infty(\Omega)},
\]
so that \eqref{eq:d_lower_bound} holds with
$b=\|a\|_{L^\infty(\Omega)}$.

An application-oriented example is provided by Arrhenius-type
reaction terms arising in combustion models; see
\cite[Eq.~(2)]{BeckerBraackVexler2004}.
For fixed parameters $A,E>0$ and $0<c<\vartheta$, consider
\[
r(t):=A\,t(c-t)\exp\left(-\frac{E}{\vartheta-t}\right),
\qquad 0\le t\le c.
\]
We retain this reaction law on the physically relevant
interval $[0,c]$ and extend it outside this interval
by its tangent lines:
\[
d(x,t):=
\begin{cases}
	r'(0)t, & t<0,\\
	r(t), & 0\le t\le c,\\
	r'(c)(t-c), & t>c.
\end{cases}
\]
Here, we have used $r(0)=r(c)=0$. The extended function belongs to $C^1(\mathbb R)$
and has a globally bounded derivative. Moreover,
\[
\partial_t d(x,t)\ge-Ac
\qquad\text{for all }t\in\mathbb R.
\]
Thus, Assumption~\ref{ass:d} holds with $b=Ac$. The modification leaves the reaction law unchanged
on $[0,c]$ and provides a globally defined nonlinearity
covered by our first-order analysis.

Under Assumption~\ref{ass:d}, both existence and uniqueness
of a solution to the state equation may fail for a given
control. To illustrate this, let
$0<\lambda_1<\lambda_2\le\lambda_3\le\cdots$
denote the eigenvalues of the Dirichlet Laplacian on
$\Omega$, counted with multiplicities, and consider
\[
d(x,t)=-\lambda_1 t.
\]
This nonlinearity satisfies Assumption~\ref{ass:d}
with $b=\lambda_1$.
For $q=0$, every multiple of a corresponding eigenfunction
$\varphi_1$ solves the state equation.
On the other hand, solvability for a given $q$ requires
the compatibility condition
\[
\int_\omega q\varphi_1\,dx=0.
\]

Nonuniqueness also occurs for nonlinear state equations.
For example, consider
\[
d(x,t)=\mu(t^3-t),
\qquad \mu>0.
\]
For $q=0$, the state equation becomes
\[
-\Delta u=\mu(u-u^3)
\quad\text{in }\Omega,
\qquad
u=0\quad\text{on }\partial\Omega.
\]
If $\mu>\lambda_1$, this equation admits at least three
distinct solutions: the zero solution, a positive solution,
and its negative.
More generally, if $\mu>\lambda_k$, it admits at least
$k$ distinct pairs of nonzero solutions, in addition to
the zero solution.
This follows by applying Clark's theorem
\cite{Clark:1972} to the associated energy functional;
see also \cite[Theorem~A]{LiWang:2002}.
For a broader discussion of existence, multiplicity,
and bifurcation of positive solutions to semilinear
elliptic equations, we refer to \cite{Lions:1982}.

Consequently, the state equation does not, in general,
define a single-valued control-to-state mapping,
and we formulate the optimization problem in terms
of admissible control-state pairs.

A further difficulty concerns the linearized state equation.
Given a state $u$, its formal linearization in a control direction
$h\in L^2(\omega)$ reads
\[
\begin{aligned}
	-\Delta v+\partial_t d(x,u)v &= h\chi_\omega
	&& \text{in }\Omega,\\
	v &=0
	&& \text{on }\partial\Omega.
\end{aligned}
\]
The lower bound \eqref{eq:d_lower_bound} does not ensure coercivity
or invertibility of the linearized operator. In particular, this
operator may have a nontrivial kernel, in which case solutions
of the linearized equation need not exist and, if they exist,
need not be unique. This phenomenon also occurs for genuinely
nonlinear equations mentioned above: for
\[
d(x,t)=\mu(t^3-t),
\]
the linearization at the solution $u=0$, corresponding to $q=0$,
is $-\Delta-\mu$, which is singular whenever $\mu$ is a Dirichlet
eigenvalue of $-\Delta$.

The same difficulty arises for the adjoint equation. Thus,
the derivation of an optimality system requires an argument
that does not rely on the invertibility of the linearized
state operator.

We denote the usual scalar products on $L^2(\Omega)$ and
$L^2(\omega)$ by $(\cdot,\cdot)$ and $(\cdot,\cdot)_\omega$,
respectively, and use the same notation for vector-valued functions. Whenever the product $fg$ belongs to $L^1(\Omega)$,
we also use the notation
\[
(f,g):=\int_\Omega fg\,dx.
\]
For a given $u\in L^1(\Omega)$, we write
\[
d(u):=d(\cdot,u(\cdot)).
\]
Since $d$ satisfies the Carath\'eodory conditions, the function
$d(u)$ is measurable; see \cite[Lemma~1]{AppellZabrejko:1989}. Since no growth conditions are imposed on $d$, the function
$d(u)$ need not belong to $L^2(\Omega)$.
We therefore also use the notation
\[
(d(u),\varphi):=\int_\Omega d(x,u(x))\varphi(x)\,dx
\]
whenever $d(u)\varphi\in L^1(\Omega)$.

For later use, we introduce the modified nonlinearity
\[
\widetilde d(x,t):=d(x,t)+bt.
\]
By \eqref{eq:d_lower_bound}, the mapping
$t\mapsto\widetilde d(x,t)$ is monotonically nondecreasing
for almost every $x\in\Omega$.
We consider the auxiliary equation
\begin{equation}\label{eq:monotone_equation}
	\begin{aligned}
		-\Delta w+\widetilde d(x,w) &= f
		&& \text{in }\Omega,\\
		w &=0
		&& \text{on }\partial\Omega.
	\end{aligned}
\end{equation}
For a given $f \in L^p(\Omega)$ with $p>N/2$, a function $w \in H^1_0(\Omega)$ is called a weak solution of \eqref{eq:monotone_equation} if $\tilde d(w)\varphi \in L^1(\Omega)$ for every $\varphi\in H^1_0(\Omega)$
and
\begin{equation}\label{eq:weak_state_modified}
	(\nabla w,\nabla\varphi)+(\tilde d(w),\varphi)
	=(f,\varphi)
	\qquad\text{for all }\varphi\in H^1_0(\Omega).
\end{equation}
The following result follows, e.g., from the arguments in
\cite[Section~2.6, in particular Theorem~2.6.11]{VexlerMeidner:2025},
which extend to the present assumption
$d(\cdot,t)\in L^2(\Omega)$.
For the Stampacchia truncation argument underlying the
$L^\infty(\Omega)$ estimate, we refer to
\cite[Appendix~B]{KinderlehrerStampacchia:2000};
continuity up to the boundary follows from the regularity
theory for linear elliptic equations, see \cite{Stampacchia:1965}.

\begin{lemma}\label{lem:monotone_equation}
	Let $p>N/2$ and $f\in L^p(\Omega)$.
	Then \eqref{eq:monotone_equation} admits a unique weak solution
	$w\in H^1_0(\Omega)\cap C(\overline\Omega)$.
	Moreover,
	\begin{equation}\label{eq:monotone_estimate}
		\|w\|_{H^1_0(\Omega)}
		+\|w\|_{C(\overline\Omega)}
		\le
		c\left(
		\|f\|_{L^p(\Omega)}
		+\|d(\cdot,0)\|_{L^2(\Omega)}
		\right),
	\end{equation}
	where $c$ depends only on $\Omega$, $N$, and $p$,
	but is independent of $d$, $b$, and $f$.
\end{lemma}

We will also use the following local Lipschitz estimates.
\begin{lemma}\label{lemma:d_Lipschitz}
	Let Assumption \ref{ass:d} be fulfilled. Let $M>0$, $u_1,u_2 \in L^\infty(\Omega)$ with $\|u_i\|_{L^\infty(\Omega)} \le M$ for $i=1,2$. Then there exists a constant $C_M>0$, independent of
	$u_1$ and $u_2$, such that
	\[
	\|d(u_1)-d(u_2)\|_{L^2(\Omega)} \le C_M \|u_1-u_2\|_{L^2(\Omega)}.
	\]
	If in addition Assumption \ref{ass:d_second_order} holds, then there holds
	\[
	\|\partial_t d(u_1)-\partial_t d(u_2)\|_{L^2(\Omega)} \le C_M \|u_1-u_2\|_{L^2(\Omega)}
	\]
	and
	\[
	\|\partial_{tt} d(u_1)-\partial_{tt} d(u_2)\|_{L^2(\Omega)} \le C_M \|u_1-u_2\|_{L^2(\Omega)}
	\]
\end{lemma}
\begin{proof}
The proof follows from Assumption \ref{ass:d} and Assumption \ref{ass:d_second_order} by a direct calculation.
\end{proof}
\section{Existence of an optimal pair}\label{sec:existence}

In this section, we prove the existence of a globally optimal
control-state pair.

Analogously to \eqref{eq:weak_state_modified}, we define
weak solutions of the state equation \eqref{eq:ocp_state}
as follows. For a given $q\in L^2(\omega)$, a function
$u\in H^1_0(\Omega)$ is called a weak solution of
\eqref{eq:ocp_state} if
$d(u)\varphi\in L^1(\Omega)$ for every $\varphi\in H^1_0(\Omega)$
and
\begin{equation}\label{eq:weak_state}
	(\nabla u,\nabla\varphi)+(d(u),\varphi)
	=(q,\varphi)_\omega
	\qquad\text{for all }\varphi\in H^1_0(\Omega).
\end{equation}
We define the set of admissible control-state pairs by
\[
W_{\mathrm{ad}}
:=
\left\{
(q,u)\in Q_{\mathrm{ad}}\times H^1_0(\Omega):
u \text{ is a weak solution of \eqref{eq:ocp_state}}
\right\}.
\]

\begin{lemma}\label{lem:admissible_state_estimate}
	There exists a constant $C>0$, depending only on $\Omega$,
	such that every admissible pair $(q,u)\in W_{\mathrm{ad}}$
	satisfies $u\in C(\overline\Omega)$ and
	\[
	\|\nabla u\|_{L^2(\Omega)}
	+\|u\|_{L^\infty(\Omega)}
	\le C\left(
	\|q\|_{L^2(\omega)}
	+b\|u\|_{L^2(\Omega)}
	+\|d(0)\|_{L^2(\Omega)}
	\right).
	\]
\end{lemma}
\begin{proof}
	Since $(q,u)\in W_{\mathrm{ad}}$, the state $u$ is a weak
	solution of \eqref{eq:weak_state}. Thus, it is also a weak solution of \eqref{eq:monotone_equation} with
	$f=q\chi_\omega+bu\in L^2(\Omega)$.
	The assertion therefore follows from
	Lemma~\ref{lem:monotone_equation}.
\end{proof}

\begin{theorem}\label{thm:existence}
Let the standing assumptions on the problem data and
Assumption~\ref{ass:d} be satisfied, and suppose that
$W_{\mathrm{ad}}\neq\emptyset$.
	Then there exists a globally optimal pair
	$(\bar q,\bar u)\in W_{\mathrm{ad}}$, that is,
	\[
	J(\bar q,\bar u)
	=
	\min_{(q,u)\in W_{\mathrm{ad}}}J(q,u).
	\]
	Moreover, every globally optimal pair satisfies
	$\bar u\in H^1_0(\Omega)\cap C(\overline\Omega)$.
\end{theorem}
\begin{proof}
	Since $W_{\mathrm{ad}}\neq\emptyset$ and $J$ is nonnegative,
	the value
	\[
	\bar J:=\inf_{(q,u)\in W_{\mathrm{ad}}}J(q,u)
	\]
	is finite. Let $\{(q_n,u_n)\}\subset W_{\mathrm{ad}}$
	be a minimizing sequence, that is,
	\[
	J(q_n,u_n)\to \bar J.
	\]
	The definition of $J$ and $\alpha>0$ imply that
	$\{q_n\}$ is bounded in $L^2(\omega)$ and
	$\{u_n\}$ is bounded in $L^2(\Omega)$.
	Consequently, after passing to a subsequence, we have
	\[
	q_n\rightharpoonup\bar q
	\quad\text{in }L^2(\omega)
	\]
	for some $\bar q\in L^2(\omega)$.
	Since $Q_{\mathrm{ad}}$ is closed and convex in $L^2(\omega)$,
	it is weakly closed, and hence $\bar q\in Q_{\mathrm{ad}}$. From $(q_n,u_n)\in W_{\mathrm{ad}}$ we obtain by \cref{lem:admissible_state_estimate}	
	\[
	u_n\in H^1_0(\Omega)\cap C(\overline\Omega)
	\]
	and
	\[
	\|u_n\|_{H^1_0(\Omega)}
	+\|u_n\|_{C(\overline\Omega)}
	\le
	c\left(
	\|q_n\|_{L^2(\omega)}
	+b\|u_n\|_{L^2(\Omega)}
	+\|d(0)\|_{L^2(\Omega)}
	\right).
	\]
	Thus, $\{u_n\}$ is bounded in
	$H^1_0(\Omega)\cap L^\infty(\Omega)$.
	By compactness, we may extract a further subsequence such that
	\[
	\begin{aligned}
		u_n&\rightharpoonup\bar u
		&&\text{in }H^1_0(\Omega),\\
		u_n&\to\bar u
		&&\text{in }L^2(\Omega)
		\text{ and almost everywhere in }\Omega.
	\end{aligned}
	\]
	Moreover, there exists $M>0$ such that
	$|u_n|\le M$ and $|\bar u|\le M$ almost everywhere in $\Omega$. Thus, \cref{lemma:d_Lipschitz} implies
	\[
	d(u_n)\to d(\bar u)
	\qquad\text{in }L^2(\Omega).
	\]
	Passing to the limit in \eqref{eq:weak_state}, we obtain
	\[
	(\nabla\bar u,\nabla\varphi)+(d(\bar u),\varphi)
	=(\bar q,\varphi)_\omega
	\qquad\text{for all }\varphi\in H^1_0(\Omega).
	\]
	Therefore, $(\bar q,\bar u)\in W_{\mathrm{ad}}$.
	By weak lower semicontinuity,
	\[
	J(\bar q,\bar u)
	\le\liminf_{n\to\infty}J(q_n,u_n)
	=\bar J,
	\]
	which proves global optimality. Finally, \cref{lem:admissible_state_estimate} implies that
	the state of every admissible pair, and hence of every
	globally optimal pair, belongs to
	$H^1_0(\Omega)\cap C(\overline\Omega)$.
\end{proof}
\begin{remark}
	The assumption $W_{\mathrm{ad}}\neq\emptyset$ cannot be
	omitted in general. Indeed, let $d(x,t)=-\lambda_1t$
	and $q_a=1$, $q_b=2$. Testing the state equation with
	a positive first Dirichlet eigenfunction $\varphi_1$
	would yield
	\[
	0=(q,\varphi_1)_\omega
	\ge\int_\omega\varphi_1\,dx>0,
	\]
	which is impossible.
	On the other hand, if $d(\cdot,0)=0$ and
	$q_a\le0\le q_b$, then $(0,0)\in W_{\mathrm{ad}}$.
\end{remark}

Besides globally optimal pairs, we will also consider local solutions
with respect to the $L^2(\omega)\times L^2(\Omega)$ topology.

\begin{definition}\label{def:local_solution}
	A pair $(\bar q,\bar u)\in W_{\mathrm{ad}}$ is called a local
	solution of \eqref{eq:ocp} if there exists $\varepsilon>0$ such that
	\[
	J(\bar q,\bar u)\le J(q,u)
	\]
	for all $(q,u)\in W_{\mathrm{ad}}$ satisfying
	\[
	\|q-\bar q\|_{L^2(\omega)}
	+\|u-\bar u\|_{L^2(\Omega)}<\varepsilon.
	\]
\end{definition}

\section{Penalized problem}\label{sec:penalized}

In this section, we introduce a family of penalized problems
whose state equations are well posed.
These problems will be used to derive first-order optimality
conditions for a given local solution $(\bar q,\bar u)$
of \eqref{eq:ocp}.

Using the modified nonlinearity $\widetilde d(x,t)=d(x,t)+bt$,
we first rewrite \eqref{eq:ocp} equivalently as
\begin{subequations}\label{eq:reformulated_ocp}
	\begin{equation}
		\text{Minimize}\quad J(q,u)
		\quad\text{over }(q,r,u)\in
		Q_{\mathrm{ad}}\times L^2(\Omega)\times H^1_0(\Omega),
	\end{equation}
	subject to
	\begin{equation}\label{eq:auxiliary_state}
		\begin{aligned}
			-\Delta u+\widetilde d(x,u)&=q\chi_\omega+r
			&&\text{in }\Omega,\\
			u&=0
			&&\text{on }\partial\Omega,
		\end{aligned}
	\end{equation}
	and the additional constraint
	\begin{equation}\label{eq:mixed_constraint}
		bu-r=0
		\qquad\text{in }L^2(\Omega).
	\end{equation}
\end{subequations}
The state equation is understood in the weak sense.
The equivalence follows by associating each
$(q,u)\in W_{\mathrm{ad}}$ with the triple $(q,bu,u)$.
In contrast to the original state equation,
\eqref{eq:auxiliary_state} admits a unique solution
$u\in H^1_0(\Omega)\cap C(\overline\Omega)$
for every $(q,r)\in Q_{\mathrm{ad}}\times L^2(\Omega)$
by Lemma~\ref{lem:monotone_equation}.

We now penalize \eqref{eq:mixed_constraint}.
Fix a local solution $(\bar q,\bar u)$ of \eqref{eq:ocp}.
For $\gamma>0$ and $\sigma\ge0$, define
\begin{equation}\label{eq:penalized_functional}
	\begin{aligned}
		J_{\gamma,\sigma}(q,r,u)
		:={}&J(q,u)
		+\frac{\gamma}{2}\|bu-r\|_{L^2(\Omega)}^2\\
		&+\frac{\sigma}{2}\|q-\bar q\|_{L^2(\omega)}^2
		+\frac{\sigma}{2}\|u-\bar u\|_{L^2(\Omega)}^2.
	\end{aligned}
\end{equation}
The corresponding penalized problem is
\begin{equation}\label{eq:penalized_problem}
	\begin{aligned}
		&\text{Minimize}\quad J_{\gamma,\sigma}(q,r,u)\\
		&\text{over }(q,r,u)\in
		Q_{\mathrm{ad}}\times L^2(\Omega)\times H^1_0(\Omega)
		\quad\text{subject to }\eqref{eq:auxiliary_state}.
	\end{aligned}
\end{equation}
The additional terms for $\sigma>0$ implement Barbu's
localization technique: they will allow us to approximate the
prescribed local solution $(\bar q,\bar u)$ by global solutions
of the penalized problems when $\sigma$ is chosen sufficiently large.

\begin{lemma}\label{lem:penalized_existence}
	For every $\gamma>0$ and $\sigma\ge0$, the penalized problem
	\eqref{eq:penalized_problem} admits a global solution
	\[
	(q_\gamma,r_\gamma,u_\gamma)
	\in Q_{\mathrm{ad}}\times L^2(\Omega)\times
	\bigl(H^1_0(\Omega)\cap C(\overline\Omega)\bigr).
	\]
	Here and below, we suppress the dependence on $\sigma$ in the notation.
\end{lemma}

\begin{proof}
	The assertion follows by standard arguments of the direct method
	of the calculus of variations, using the monotonicity of
	$\widetilde d$ and Lemma~\ref{lem:monotone_equation}
	for the existence, uniqueness, and a priori bounds
	of solutions to \eqref{eq:auxiliary_state}.
\end{proof}

\begin{lemma}\label{lem:penalized_convergence}
	Let $(\bar q,\bar u)$ be a local solution of \eqref{eq:ocp}
	in the sense of Definition~\ref{def:local_solution},
	and let $\varepsilon>0$ be a corresponding radius.
	Choose $\sigma>0$ such that
	\[
	\sigma>\frac{4J(\bar q,\bar u)}{\varepsilon^2}.
	\]
	For each $\gamma>0$, let $(q_\gamma,r_\gamma,u_\gamma)$
	be any global solution of \eqref{eq:penalized_problem}.
	Then, as $\gamma\to\infty$,
	\[
	q_\gamma\to\bar q
	\quad\text{strongly in }L^2(\omega),
	\qquad
	u_\gamma\to\bar u
	\quad\text{strongly in }L^2(\Omega).
	\]
	Moreover,
	\[
	r_\gamma\to b\bar u
	\qquad\text{strongly in }L^2(\Omega).
	\]
\end{lemma}
\begin{proof}
	Since $(\bar q,b\bar u,\bar u)$ is admissible for the penalized
	problem, global optimality gives
	\begin{equation}\label{eq:penalty_comparison}
		\begin{aligned}
			J(q_\gamma,u_\gamma)
			&+\frac{\gamma}{2}\|bu_\gamma-r_\gamma\|_{L^2(\Omega)}^2\\
			&+\frac{\sigma}{2}\|q_\gamma-\bar q\|_{L^2(\omega)}^2
			+\frac{\sigma}{2}\|u_\gamma-\bar u\|_{L^2(\Omega)}^2
			\le J(\bar q,\bar u).
		\end{aligned}
	\end{equation}
	In particular, $\{q_\gamma\}$ and $\{u_\gamma\}$ are bounded
	in $L^2(\omega)$ and $L^2(\Omega)$, respectively, and
	\[
	\|bu_\gamma-r_\gamma\|_{L^2(\Omega)}
	\le \left(\frac{2J(\bar q,\bar u)}{\gamma}\right)^{1/2}
	\to0.
	\]
	Consequently, $\{r_\gamma\}$ is bounded in $L^2(\Omega)$
	as $\gamma\to\infty$.
	Applying Lemma~\ref{lem:monotone_equation} to
	\eqref{eq:auxiliary_state} shows that $\{u_\gamma\}$ is
	uniformly bounded in $H^1_0(\Omega)\cap C(\overline\Omega)$.
	
	Let $\gamma_n\to\infty$ be arbitrary.
	After passing to a subsequence, we obtain
	\[
	\begin{aligned}
		q_{\gamma_n}&\rightharpoonup q^*
		&&\text{in }L^2(\omega),\\
		u_{\gamma_n}&\rightharpoonup u^*
		&&\text{in }H^1_0(\Omega),\\
		u_{\gamma_n}&\to u^*
		&&\text{in }L^2(\Omega)
		\text{ and almost everywhere in }\Omega.
	\end{aligned}
	\]
	Moreover, $r_{\gamma_n}\to bu^*$ in $L^2(\Omega)$.
	As in the proof of Theorem~\ref{thm:existence}, \cref{lemma:d_Lipschitz} as well as 
	the uniform boundedness of $u_{\gamma_n}$ and $u^*$
	imply
	\[
	\widetilde d(u_{\gamma_n})\to\widetilde d(u^*)
	\qquad\text{in }L^2(\Omega).
	\]
	Passing to the limit in the auxiliary state equation gives
	\[
	(\nabla u^*,\nabla\varphi)+(\widetilde d(u^*),\varphi)
	=(q^*,\varphi)_\omega+(bu^*,\varphi)
	\qquad\text{for all }\varphi\in H^1_0(\Omega).
	\]
	Since $Q_{\mathrm{ad}}$ is weakly closed, this implies
	$(q^*,u^*)\in W_{\mathrm{ad}}$.
	
	By \eqref{eq:penalty_comparison} and weak lower
	semicontinuity,
	\[
	\|q^*-\bar q\|_{L^2(\omega)}
	+\|u^*-\bar u\|_{L^2(\Omega)}
	\le
	2\sqrt{\frac{J(\bar q,\bar u)}{\sigma}}
	<\varepsilon.
	\]
	Local optimality of $(\bar q,\bar u)$ therefore yields
	\[
	J(\bar q,\bar u)\le J(q^*,u^*).
	\]
On the other hand, weak lower semicontinuity and
\eqref{eq:penalty_comparison} yield
\[
\begin{aligned}
	J(q^*,u^*)
	&+\frac{\sigma}{2}\|q^*-\bar q\|_{L^2(\omega)}^2
	+\frac{\sigma}{2}\|u^*-\bar u\|_{L^2(\Omega)}^2\\
	&\le \liminf_{n\to\infty}
	\Bigl(
	J(q_{\gamma_n},u_{\gamma_n})
	+\frac{\sigma}{2}\|q_{\gamma_n}-\bar q\|_{L^2(\omega)}^2
	+\frac{\sigma}{2}\|u_{\gamma_n}-\bar u\|_{L^2(\Omega)}^2
	\Bigr)\\
	&\le J(\bar q,\bar u).
\end{aligned}
\]
	Since $\sigma>0$, we conclude that
	$(q^*,u^*)=(\bar q,\bar u)$.
	
We already know that $u_{\gamma_n}\to\bar u$ strongly in
$L^2(\Omega)$ and $q_{\gamma_n}\rightharpoonup\bar q$
weakly in $L^2(\omega)$.
To prove strong convergence of the controls, we use
\eqref{eq:penalty_comparison} to obtain
\[
\frac{\sigma}{2}
\|q_{\gamma_n}-\bar q\|_{L^2(\omega)}^2
\le J(\bar q,\bar u)-J(q_{\gamma_n},u_{\gamma_n}).
\]
Consequently, weak lower semicontinuity of $J$ yields
\[
\begin{aligned}
	0
	&\le \frac{\sigma}{2}\limsup_{n\to\infty}
	\|q_{\gamma_n}-\bar q\|_{L^2(\omega)}^2\\
	&\le J(\bar q,\bar u)
	-\liminf_{n\to\infty}J(q_{\gamma_n},u_{\gamma_n})
	\le 0.
\end{aligned}
\]
Since $\sigma>0$, this proves
$q_{\gamma_n}\to\bar q$ strongly in $L^2(\omega)$.

Thus, every sequence $\gamma_n\to\infty$ admits a subsequence
along which $(q_{\gamma_n},u_{\gamma_n})$ converges strongly
to $(\bar q,\bar u)$ in $L^2(\omega)\times L^2(\Omega)$.
This implies the asserted convergence as $\gamma\to\infty$. The convergence of $r_\gamma$ follows from
$r_\gamma-b\bar u
=b(u_\gamma-\bar u)-(bu_\gamma-r_\gamma)$.
\end{proof}

We next state the first-order optimality conditions for the
penalized problem.

\begin{lemma}\label{lem:penalized_optimality}
	Let $\gamma>0$ and $\sigma\ge0$, and let
	$(q_\gamma,r_\gamma,u_\gamma)$ be a global solution of
	\eqref{eq:penalized_problem}.
	Then there exists a unique adjoint state
	$z_\gamma\in H^1_0(\Omega)\cap C(\overline\Omega)$
	such that the following optimality system holds:
	\begin{subequations}\label{eq:penalized_optimality_system}
		\begin{equation}\label{eq:penalized_state}
			(\nabla u_\gamma,\nabla\varphi)
			+(\widetilde d(u_\gamma),\varphi)
			=(q_\gamma,\varphi)_\omega+(r_\gamma,\varphi)
		\end{equation}
		and
		\begin{equation}\label{eq:penalized_adjoint}
			\begin{aligned}
				(\nabla z_\gamma,\nabla\varphi)
				+(\partial_t\widetilde d(\cdot,u_\gamma)z_\gamma,\varphi)
				={}&(u_\gamma-u_d+\sigma(u_\gamma-\bar u),\varphi)\\
				&+\gamma b(bu_\gamma-r_\gamma,\varphi)
			\end{aligned}
		\end{equation}
		for all $\varphi\in H^1_0(\Omega)$, together with
		$q_\gamma\in Q_{\mathrm{ad}}$,
		\begin{equation}\label{eq:penalized_control_vi}
			\bigl(
			\alpha q_\gamma+\sigma(q_\gamma-\bar q)+z_\gamma|_\omega,\,
			q-q_\gamma
			\bigr)_\omega\ge0
			\qquad\text{for all }q\in Q_{\mathrm{ad}},
		\end{equation}
		and
		\begin{equation}\label{eq:penalized_r_optimality}
			z_\gamma=\gamma(bu_\gamma-r_\gamma)
			\qquad\text{in }L^2(\Omega).
		\end{equation}
	\end{subequations}
\end{lemma}

\begin{remark}
Using \eqref{eq:penalized_r_optimality} and
$\partial_t\widetilde d=\partial_t d+b$,
the adjoint equation \eqref{eq:penalized_adjoint}
can equivalently be written as
\begin{equation}\label{eq:penalized_adjoint_simplified}
	(\nabla z_\gamma,\nabla\varphi)
	+(\partial_t d(\cdot,u_\gamma)z_\gamma,\varphi)
	=
	(u_\gamma-u_d+\sigma(u_\gamma-\bar u),\varphi)
\end{equation}
for all $\varphi\in H^1_0(\Omega)$.
\end{remark}

\begin{proof}[Proof of \cref{lem:penalized_optimality}]
	The result follows by standard arguments for optimal control
	problems governed by monotone semilinear elliptic equations;
	see \cite[Section~7.2, in particular Theorem~7.2.11]{VexlerMeidner:2025}.
	The arguments apply with the control pair
	$(q,r)\in Q_{\mathrm{ad}}\times L^2(\Omega)$
	and the modified cost functional $J_{\gamma,\sigma}$.
	Variations with respect to $q$ yield
	\eqref{eq:penalized_control_vi}, whereas unrestricted
	variations with respect to $r$ give
	\eqref{eq:penalized_r_optimality}.
\end{proof}

\section{Optimality system}\label{sec:optimality_system}

In this section, we derive first-order optimality conditions
for a local solution $(\bar q,\bar u)$ of \eqref{eq:ocp}
in the sense of Definition~\ref{def:local_solution}.
We define the inactive set associated with $\bar q$ by
\[
I(\bar q):=\{x\in\omega:q_a<\bar q(x)<q_b\}.
\]
This set is understood up to sets of Lebesgue measure zero.
We impose the following assumption.

\begin{assumption}\label{ass:inactive_set}
	The inactive set has positive Lebesgue measure, that is,
	\[
	|I(\bar q)|>0.
	\]
\end{assumption}

In particular, Assumption~\ref{ass:inactive_set} is automatically
satisfied in the unconstrained case $q_a=-\infty$ and $q_b=+\infty$,
since then $I(\bar q)=\omega$ up to a set of measure zero.

The derivation of the optimality system relies on the following
unique continuation principle for elliptic equations with
bounded potentials.
In particular, vanishing on a measurable set of positive
Lebesgue measure is sufficient.

\begin{lemma}[Unique continuation]\label{lem:unique_continuation}
	Let $\Omega\subset\mathbb{R}^N$ be open and connected,
	and let $a\in L^\infty(\Omega)$.
	Suppose that $z\in H^1_{\mathrm{loc}}(\Omega)$ satisfies
	\[
	-\Delta z+az=0
	\qquad\text{in }\Omega
	\]
	in the sense of distributions.
	If there exists a measurable set $E\subset\Omega$ with
	$|E|>0$ such that $z=0$ almost everywhere on $E$,
	then $z=0$ almost everywhere in $\Omega$.
\end{lemma}
\begin{proof}
	The assertion follows, e.g., from
	\cite[Theorem~2.4(a)]{HarrachPohjolaSalo:2019}
	by localization and connectedness of $\Omega$.
	For a classical treatment in dimensions $N\ge3$,
	we also refer to
	\cite[Proposition~3 and Remark~2]{deFigueiredoGossez:1992}.
\end{proof}

We are now in a position to state the first of our main results.
Under Assumption~\ref{ass:inactive_set}, we obtain a normal
optimality system without requiring invertibility of the
linearized operator $-\Delta+\partial_t d(\cdot,\bar u)$.

\begin{theorem}\label{thm:optimality_system}
	Let the standing assumptions on the problem data and
	Assumption~\ref{ass:d} be satisfied.
	Let $(\bar q,\bar u)\in W_{\mathrm{ad}}$ be a local solution
	of \eqref{eq:ocp} in the sense of
	Definition~\ref{def:local_solution}.
	Suppose that Assumption~\ref{ass:inactive_set} holds.
	Then there exists a unique adjoint state
	$\bar z\in H^1_0(\Omega)\cap C(\overline\Omega)$
	such that the following optimality system holds:
	\begin{subequations}\label{eq:optimality_system}
		\begin{align}
			(\nabla\bar u,\nabla\varphi)+(d(\bar u),\varphi)
			&=(\bar q,\varphi)_\omega,
			\label{eq:optimal_state}\\
			(\nabla\bar z,\nabla\varphi)
			+(\partial_t d(\cdot,\bar u)\bar z,\varphi)
			&=(\bar u-u_d,\varphi),
			\label{eq:optimal_adjoint}
		\end{align}
		for all $\varphi\in H^1_0(\Omega)$, together with
		\begin{equation}\label{eq:optimal_control}
			\bar q\in Q_{\mathrm{ad}},
			\qquad
			(\alpha\bar q+\bar z|_\omega,q-\bar q)_\omega\ge0
			\quad\text{for all }q\in Q_{\mathrm{ad}}.
		\end{equation}
	\end{subequations}
	Equivalently, the control condition can be expressed as
	\begin{equation}\label{eq:optimal_projection}
		\bar q(x)
		=
		\operatorname{proj}_{[q_a,q_b]}
		\left(-\frac{\bar z(x)}{\alpha}\right)
		\qquad\text{for almost every }x\in\omega.
	\end{equation}
\end{theorem}
\begin{proof}
	\textit{Step 1: Convergence of the penalized solutions and
		preliminary estimates.}
	Fix $\sigma>0$ sufficiently large as in
	Lemma~\ref{lem:penalized_convergence}.
	For each $\gamma>0$, let $(q_\gamma,r_\gamma,u_\gamma)$
	be a global solution of \eqref{eq:penalized_problem},
	and let $z_\gamma$ be the corresponding adjoint state
	from Lemma~\ref{lem:penalized_optimality}.
	By Lemma~\ref{lem:penalized_convergence},
	\[
	q_\gamma\to\bar q
	\quad\text{in }L^2(\omega),
	\qquad
	u_\gamma\to\bar u
	\quad\text{in }L^2(\Omega).
	\]
	Moreover, the estimates in its proof yield a constant $M>0$
	such that, for all sufficiently large $\gamma$,
	\[
	\|u_\gamma\|_{L^\infty(\Omega)}\le M,
	\qquad
	\|\bar u\|_{L^\infty(\Omega)}\le M.
	\]
	Define
	\[
	a_\gamma:=\partial_t d(\cdot,u_\gamma),
	\qquad
	a:=\partial_t d(\cdot,\bar u).
	\]
	By Assumption~\ref{ass:d},
	\[
	\|a_\gamma\|_{L^\infty(\Omega)}\le c_M,
	\qquad
	\|a\|_{L^\infty(\Omega)}\le c_M,
	\qquad
	a_\gamma,a\ge-b
	\quad\text{almost everywhere in }\Omega.
	\]
	
	For every sequence $\gamma_n\to\infty$, we may pass to
	a subsequence such that $u_{\gamma_n}\to\bar u$
	almost everywhere in $\Omega$.
	The continuity of $\partial_t d(x,\cdot)$ for almost
	every $x\in\Omega$ then implies
	\[
	a_{\gamma_n}\to a
	\quad\text{almost everywhere in }\Omega.
	\]
	In particular, along such a subsequence, if $w_n\to w$
	in $L^2(\Omega)$, then
	\begin{equation}\label{eq:coefficient_product_convergence}
		\begin{aligned}
			\|a_{\gamma_n}w_n-aw\|_{L^2(\Omega)}
			&\le c_M\|w_n-w\|_{L^2(\Omega)}
			+\|(a_{\gamma_n}-a)w\|_{L^2(\Omega)}
			\longrightarrow0.
		\end{aligned}
	\end{equation}
	Indeed, the second term tends to zero by dominated
	convergence, since
	$|(a_{\gamma_n}-a)w|^2\le4c_M^2|w|^2$.
	
	Finally, set
	\[
	F_\gamma:=u_\gamma-u_d+\sigma(u_\gamma-\bar u).
	\]
	Then
	\[
	F_\gamma\to\bar u-u_d
	\quad\text{in }L^2(\Omega),
	\]
	and \eqref{eq:penalized_adjoint_simplified} reads
	\begin{equation}\label{eq:adjoint_gamma_rewritten}
		(\nabla z_\gamma,\nabla\varphi)
		+(a_\gamma z_\gamma,\varphi)
		=(F_\gamma,\varphi)
		\qquad\text{for all }\varphi\in H^1_0(\Omega).
	\end{equation}
	
	\medskip
	\textit{Step 2: Boundedness of the adjoint states.}
	We prove that $\{z_\gamma\}$ remains bounded in
	$L^2(\Omega)$ as $\gamma\to\infty$.
	Suppose, to the contrary, that there exists a sequence
	$\gamma_n\to\infty$ such that
	\[
	t_n:=\|z_{\gamma_n}\|_{L^2(\Omega)}\to\infty.
	\]
	Passing to a subsequence as in Step~1, we may assume
	that $a_{\gamma_n}\to a$ almost everywhere in $\Omega$.
	Define
	\[
	v_n:=\frac{z_{\gamma_n}}{t_n},
	\qquad
	\|v_n\|_{L^2(\Omega)}=1.
	\]
	Dividing \eqref{eq:adjoint_gamma_rewritten} by $t_n$ gives
	\[
	(\nabla v_n,\nabla\varphi)
	+(a_{\gamma_n}v_n,\varphi)
	=t_n^{-1}(F_{\gamma_n},\varphi)
	\qquad\text{for all }\varphi\in H^1_0(\Omega).
	\]
	Testing with $\varphi=v_n$ and using $a_{\gamma_n}\ge-b$,
	we obtain
	\[
	\|\nabla v_n\|_{L^2(\Omega)}^2
	\le b+t_n^{-1}\|F_{\gamma_n}\|_{L^2(\Omega)}.
	\]
	Thus, $\{v_n\}$ is bounded in $H^1_0(\Omega)$.
	After passing to a further subsequence,
	\[
	v_n\rightharpoonup v
	\quad\text{in }H^1_0(\Omega),
	\qquad
	v_n\to v
	\quad\text{in }L^2(\Omega),
	\qquad
	\|v\|_{L^2(\Omega)}=1.
	\]
	By \eqref{eq:coefficient_product_convergence},
	$a_{\gamma_n}v_n\to av$ in $L^2(\Omega)$.
	Passing to the limit in the normalized adjoint equation
	therefore yields
	\begin{equation}\label{eq:normalized_adjoint_limit}
		(\nabla v,\nabla\varphi)+(av,\varphi)=0
		\qquad\text{for all }\varphi\in H^1_0(\Omega).
	\end{equation}
	
	Next, dividing \eqref{eq:penalized_control_vi} by $t_n$,
	we obtain
	\[
	\left(
	t_n^{-1}\bigl(\alpha q_{\gamma_n}
	+\sigma(q_{\gamma_n}-\bar q)\bigr)
	+v_n|_\omega,\,
	q-q_{\gamma_n}
	\right)_\omega\ge0
	\qquad\text{for all }q\in Q_{\mathrm{ad}}.
	\]
	Using the strong convergence of $q_{\gamma_n}$ and $v_n$,
	we pass to the limit and conclude that
	\begin{equation}\label{eq:normalized_control_vi}
		(v,q-\bar q)_\omega\ge0
		\qquad\text{for all }q\in Q_{\mathrm{ad}}.
	\end{equation}
	This implies that $v=0$ almost everywhere on $I(\bar q)$.
	Indeed, define
	\[
	I_k:=
	\left\{
	x\in\omega:
	q_a+\frac1k\le\bar q(x)\le q_b-\frac1k
	\right\},
	\qquad k\in\mathbb N,
	\]
	with the usual convention for infinite bounds.
	For every measurable subset $E\subset I_k$, both controls
	$\bar q\pm k^{-1}\chi_E$ belong to $Q_{\mathrm{ad}}$.
	Inserting them into \eqref{eq:normalized_control_vi} gives
	\[
	\int_E v\,dx=0.
	\]
	Hence $v=0$ almost everywhere on each $I_k$, and therefore
	on $I(\bar q)=\bigcup_{k\in\mathbb N}I_k$ up to null sets.
	
	By Assumption~\ref{ass:inactive_set}, this zero set has
	positive Lebesgue measure.
	Since $a\in L^\infty(\Omega)$ and $\Omega$ is connected,
	Lemma~\ref{lem:unique_continuation} applied to
	\eqref{eq:normalized_adjoint_limit} implies $v=0$
	almost everywhere in $\Omega$.
	This contradicts $\|v\|_{L^2(\Omega)}=1$ and proves
	the asserted boundedness of $\{z_\gamma\}$ in $L^2(\Omega)$.
	
	\medskip
	\textit{Step 3: Passage to the limit in the optimality system.}
	Testing \eqref{eq:adjoint_gamma_rewritten} with $z_\gamma$
	gives
	\[
	\|\nabla z_\gamma\|_{L^2(\Omega)}^2
	\le b\|z_\gamma\|_{L^2(\Omega)}^2
	+\|F_\gamma\|_{L^2(\Omega)}
	\|z_\gamma\|_{L^2(\Omega)}.
	\]
	Thus, $\{z_\gamma\}$ is also bounded in $H^1_0(\Omega)$.
	Along a sequence $\gamma_n\to\infty$, we may therefore
	assume, after passing to a subsequence, that
	\[
	z_{\gamma_n}\rightharpoonup\bar z
	\quad\text{in }H^1_0(\Omega),
	\qquad
	z_{\gamma_n}\to\bar z
	\quad\text{in }L^2(\Omega),
	\]
	and $a_{\gamma_n}\to a$ almost everywhere in $\Omega$.
	By \eqref{eq:coefficient_product_convergence},
	\[
	a_{\gamma_n}z_{\gamma_n}\to a\bar z
	\quad\text{in }L^2(\Omega).
	\]
	Passing to the limit in \eqref{eq:adjoint_gamma_rewritten}
	proves \eqref{eq:optimal_adjoint}.
	Similarly, passing to the limit in
	\eqref{eq:penalized_control_vi} yields
	\[
	(\alpha\bar q+\bar z|_\omega,q-\bar q)_\omega\ge0
	\qquad\text{for all }q\in Q_{\mathrm{ad}}.
	\]
	The state equation and $\bar q\in Q_{\mathrm{ad}}$
	hold because $(\bar q,\bar u)\in W_{\mathrm{ad}}$.
	
	To obtain continuity of $\bar z$, rewrite its equation as
	\[
	-\Delta\bar z+(a+b)\bar z=\bar u-u_d+b\bar z.
	\]
	The coefficient on the left is bounded and nonnegative,
	and the right-hand side belongs to $L^2(\Omega)$.
	Lemma~\ref{lem:monotone_equation}, applied to this linear
	monotone equation, yields $\bar z\in C(\overline\Omega)$.
	
	The projection formula \eqref{eq:optimal_projection}
	follows from the characterization of the $L^2(\omega)$
	projection onto $Q_{\mathrm{ad}}$; see
	\cite[Lemma~5.2.4 and Theorem~5.2.7]{VexlerMeidner:2025}.
	
	\medskip
	\textit{Step 4: Uniqueness of the adjoint state.}
	Let $\bar z_1$ and $\bar z_2$ satisfy the optimality system,
	and set $y:=\bar z_1-\bar z_2\in H^1_0(\Omega)$.
	Subtracting the adjoint equations gives
	\[
	(\nabla y,\nabla\varphi)+(ay,\varphi)=0
	\qquad\text{for all }\varphi\in H^1_0(\Omega).
	\]
	By the projection formula,
	\[
	\bar z_1=-\alpha\bar q=\bar z_2
	\qquad\text{almost everywhere on }I(\bar q).
	\]
	Hence $y=0$ almost everywhere on $I(\bar q)$.
	Assumption~\ref{ass:inactive_set} and
	Lemma~\ref{lem:unique_continuation} imply $y=0$
	almost everywhere in $\Omega$, proving uniqueness.
\end{proof}

\begin{remark}
	Although the adjoint equation \eqref{eq:optimal_adjoint}
	need not be uniquely solvable, the adjoint state satisfying
	both \eqref{eq:optimal_adjoint} and the control variational
	inequality \eqref{eq:optimal_control} is unique under
	Assumption~\ref{ass:inactive_set}.
\end{remark}

\begin{remark}
If both control bounds are finite and $\omega$ is open
and connected, the projection formula
	with $\alpha>0$ and a continuous adjoint state excludes
	nonconstant bang-bang controls.
	In particular, for problems with a monotone nonlinearity,
	such controls cannot be optimal.
	For the nonmonotone problems considered here, however,
	nonconstant bang-bang controls may be optimal, and a normal
	optimality system may fail to exist.
	The following example illustrates this phenomenon and
	shows that Assumption~\ref{ass:inactive_set} cannot be
	omitted altogether, even when $\omega=\Omega$.
	
	Let $\Omega=\omega=(0,\pi)^2$, $q_a=-1$, $q_b=1$,
	$u_d=0$, and $\alpha>0$.
	Set
	\[
	\psi(x)=\sin(2x_1)\sin(x_2),\qquad
	q_*(x)=-\operatorname{sign}(\psi(x)),\qquad
	d(x,t)=-5t+q_*(x).
	\]
	Since $-\Delta\psi=5\psi$ and $\psi\in H^1_0(\Omega)$,
	testing the state equation with $\psi$ yields
	\[
	0=(q-q_*,\psi)
	=\int_\Omega\bigl(q\psi+|\psi|\bigr)\,dx.
	\]
	The integrand is nonnegative for every admissible control,
	and $\psi\neq0$ almost everywhere. Thus, every admissible
	pair satisfies $q=q_*$, and $(\bar q,\bar u)=(q_*,0)$
	is the unique global solution.
	The control $\bar q$ equals $-1$ on $\{x_1<\pi/2\}$
	and $1$ on $\{x_1>\pi/2\}$, so that $|I(\bar q)|=0$.
	
	A normal optimality system would require
	$(-\Delta-5)\bar z=0$ and
	\[
	\bar z\ge\alpha \quad\text{on }\{x_1<\pi/2\},
	\qquad
	\bar z\le-\alpha \quad\text{on }\{x_1>\pi/2\},
	\]
	almost everywhere.
	This contradicts the interior continuity of $\bar z$
	across $x_1=\pi/2$.
\end{remark}
\section{Second-order sufficient optimality conditions}
\label{sec:second_order}

In this section, we derive second-order sufficient optimality conditions.
We first introduce the critical cone and the quadratic form
associated with the second variation of the Lagrangian.

Throughout this section, we additionally assume that the nonlinearity 
$d(\cdot,\cdot)$ fulfills both Assumptions \ref{ass:d} and \ref{ass:d_second_order}.

In the following, let $(\bar q,\bar u)\in W_{\mathrm{ad}}$
satisfy the first-order optimality system
\eqref{eq:optimality_system} with an adjoint state $\bar z$.

The tangent cone to $Q_{\mathrm{ad}}$ at $\bar q$ is given by
\begin{equation}\label{eq:def:T_Qad}
T_{Q_{\mathrm{ad}}}(\bar q)
=
\left\{
h\in L^2(\omega):
\begin{aligned}
	h&\ge0 &&\text{a.e. on }\{\bar q=q_a\},\\
	h&\le0 &&\text{a.e. on }\{\bar q=q_b\}
\end{aligned}
\right\}.
\end{equation}
If one of the bounds is infinite, the corresponding
condition is omitted. This cone is convex and closed in $L^2(\omega)$.

\begin{definition}[Critical directions]
	\label{def:critical_cone}
	The critical cone $\mathcal C(\bar q,\bar u)$ consists of
	all pairs
	\[
	(h,v)\in L^2(\omega)\times H_0^1(\Omega)
	\]
	such that $h\in T_{Q_{\mathrm{ad}}}(\bar q)$,
	\begin{equation}\label{eq:linearized_state}
		(\nabla v,\nabla\phi)
		+(\partial_t d(\cdot,\bar u)v,\phi)
		=(h,\phi)_\omega
		\qquad\forall\phi\in H_0^1(\Omega),
	\end{equation}
	and
	\begin{equation}\label{eq:critical_direction}
		J'(\bar q,\bar u)(h,v)
		=(\bar u-u_d,v)+\alpha(\bar q,h)_\omega
		=0.
	\end{equation}
\end{definition}

By the adjoint equation and \eqref{eq:linearized_state},
the criticality condition \eqref{eq:critical_direction}
is equivalent to
\begin{equation}\label{eq:criticality_cond}
(\alpha\bar q+\bar z,h)_\omega=0.
\end{equation}

\begin{remark}
	The linearized state equation \eqref{eq:linearized_state}
	need not admit a solution for every
	$h\in T_{Q_{\mathrm{ad}}}(\bar q)$.
	Moreover, even if a solution exists, it need not be unique,
	since the operator
	$-\Delta+\partial_t d(\cdot,\bar u)$
	may have a nontrivial kernel.
	Consequently, critical directions are considered as pairs
	$(h,v)$, with the state direction $v$ included explicitly.
\end{remark}


Finally, we define the quadratic form
\begin{equation}\label{eq:second_variation}
	\begin{aligned}
		\mathcal Q(h,v)
		={}&\alpha\|h\|_{L^2(\omega)}^2
		+\|v\|_{L^2(\Omega)}^2\\
		&-\int_\Omega
		\partial_{tt}d(x,\bar u(x))\,\bar z(x)\,v(x)^2
		\,\mathrm{d}x.
	\end{aligned}
\end{equation}
This is the second variation of the Lagrangian with respect
to the control and state variables at
$(\bar q,\bar u,\bar z)$. Since $\partial_{tt}d(\cdot,\bar u)\bar z\in L^\infty(\Omega)$,
the quadratic form
\[
\mathcal Q:L^2(\omega)\times L^2(\Omega)\to\mathbb R
\]
is continuous.

The following expansion relates $\mathcal Q$ to differences
of objective values at admissible pairs.
It will be used in the proofs of both the sufficient
and the necessary second-order optimality conditions.

\begin{lemma}\label{lem:cost_expansion}
	Let Assumptions~\ref{ass:d} and~\ref{ass:d_second_order}
	hold, and let $(\bar q,\bar u)\in W_{\mathrm{ad}}$
	satisfy the first-order optimality system
	\eqref{eq:optimality_system} with an adjoint state
	$\bar z\in H^1_0(\Omega)\cap C(\overline\Omega)$.
	For every $(q,u)\in W_{\mathrm{ad}}$, set
	\[
	\delta q:=q-\bar q,
	\qquad
	\delta u:=u-\bar u.
	\]
	Then
	\begin{equation}\label{eq:cost_expansion}
		\begin{aligned}
			J(q,u)-J(\bar q,\bar u)
			={}&(\alpha\bar q+\bar z,\delta q)_\omega\\
			&+\frac12\mathcal Q(\delta q,\delta u)
			+\mathcal R(q,u).
		\end{aligned}
	\end{equation}
	For every $M>0$ such that
	\[
	\|\bar u\|_{L^\infty(\Omega)}\le M,
	\qquad
	\|u\|_{L^\infty(\Omega)}\le M,
	\]
	the remainder satisfies
	\begin{equation}\label{eq:cost_expansion_remainder}
		|\mathcal R(q,u)|
		\le
		\frac{L_M}{6}\|\bar z\|_{L^\infty(\Omega)}
		\|\delta u\|_{L^3(\Omega)}^3,
	\end{equation}
	where $L_M$ is the constant from
	Assumption~\ref{ass:d_second_order}.
\end{lemma}

\begin{proof}
	Set
	\[
	a_1(x):=\partial_t d(x,\bar u(x)),
	\qquad
	a_2(x):=\partial_{tt}d(x,\bar u(x)).
	\]
	Taylor's formula applied to $d(x,\cdot)$ gives
	\[
	d(x,u(x))-d(x,\bar u(x))
	=a_1(x)\delta u(x)
	+\frac12 a_2(x)\delta u(x)^2+\eta(x)
	\]
	for almost every $x\in\Omega$, where
	\[
	\eta(x)
	=\delta u(x)^2
	\int_0^1(1-t)
	\left[
	\partial_{tt}d\bigl(x,\bar u(x)+t\delta u(x)\bigr)
	-a_2(x)
	\right]\,dt.
	\]
	The segment between $\bar u(x)$ and $u(x)$ lies in
	$[-M,M]$.
	Hence the local uniform Lipschitz continuity of
	$\partial_{tt}d$ in
	Assumption~\ref{ass:d_second_order} implies
	\begin{equation}\label{eq:cost_expansion_pointwise_remainder}
		|\eta(x)|
		\le\frac{L_M}{6}|\delta u(x)|^3
		\qquad\text{for almost every }x\in\Omega.
	\end{equation}
	
	Testing the difference of the state equations with
	$\bar z$ yields
	\[
	\begin{aligned}
		(\nabla\delta u,\nabla\bar z)
		+(a_1\delta u,\bar z)
		+\frac12(a_2\delta u^2,\bar z)
		+(\eta,\bar z)
		=(\delta q,\bar z)_\omega.
	\end{aligned}
	\]
	On the other hand, the adjoint equation
	\eqref{eq:optimal_adjoint}, tested with $\delta u$, gives
	\[
	(\nabla\delta u,\nabla\bar z)
	+(a_1\delta u,\bar z)
	=(\bar u-u_d,\delta u).
	\]
	Combining these identities, we obtain
	\[
	(\bar u-u_d,\delta u)
	=(\bar z,\delta q)_\omega
	-\frac12(a_2\delta u^2,\bar z)-(\eta,\bar z).
	\]
	
	Since $J$ is quadratic, its exact expansion reads
	\[
	\begin{aligned}
		J(q,u)-J(\bar q,\bar u)
		={}&\alpha(\bar q,\delta q)_\omega
		+(\bar u-u_d,\delta u)\\
		&+\frac{\alpha}{2}\|\delta q\|_{L^2(\omega)}^2
		+\frac12\|\delta u\|_{L^2(\Omega)}^2\\
		={}&(\alpha\bar q+\bar z,\delta q)_\omega
		+\frac12\mathcal Q(\delta q,\delta u)
		-(\eta,\bar z).
	\end{aligned}
	\]
	Thus, \eqref{eq:cost_expansion} holds with
	\[
	\mathcal R(q,u):=-(\eta,\bar z).
	\]
	Finally, \eqref{eq:cost_expansion_pointwise_remainder}
	and $\bar z\in L^\infty(\Omega)$ imply
	\[
	|\mathcal R(q,u)|
	\le\|\bar z\|_{L^\infty(\Omega)}
	\|\eta\|_{L^1(\Omega)}
	\le
	\frac{L_M}{6}\|\bar z\|_{L^\infty(\Omega)}
	\|\delta u\|_{L^3(\Omega)}^3.
	\]
	This proves \eqref{eq:cost_expansion_remainder}.
\end{proof}

The following result shows that strict local optimality
can be established even when the state equation is not
uniquely solvable and the linearized state operator
is singular.

\begin{theorem}[Second-order sufficient optimality condition]
	\label{thm:second_order_sufficient}
		Let the standing assumptions on the problem data,
	Assumption~\ref{ass:d} and Assumption~\ref{ass:d_second_order} be satisfied. Let $(\bar q,\bar u)\in W_{\mathrm{ad}}$ fulfill the
	first-order optimality system \eqref{eq:optimality_system} with an adjoint state
	$\bar z\in H^1_0(\Omega)\cap C(\overline\Omega)$. Suppose that
	\begin{equation}\label{eq:second_order_sufficient}
		\mathcal Q(h,v)>0
		\qquad
		\forall (h,v)\in
		\mathcal C(\bar q,\bar u)\setminus\{(0,0)\}.
	\end{equation}
	Then there exist $\varepsilon>0$ and $\kappa>0$ such that
	\begin{equation}\label{eq:quadratic_growth}
		\begin{aligned}
			J(q,u)\ge J(\bar q,\bar u)
			+\frac{\kappa}{2}\Bigl(
			\|q-\bar q\|_{L^2(\omega)}^2
			+\|u-\bar u\|_{L^2(\Omega)}^2
			\Bigr)
		\end{aligned}
	\end{equation}
	for all $(q,u)\in W_{\mathrm{ad}}$ satisfying
	\[
	\|q-\bar q\|_{L^2(\omega)}
	+\|u-\bar u\|_{L^2(\Omega)}
	\le\varepsilon.
	\]
	In particular, $(\bar q,\bar u)$ is a strict local solution
	with respect to the
	$L^2(\omega)\times L^2(\Omega)$ topology.
\end{theorem}
\begin{proof}
	\textit{Step 1: Contradiction argument and preliminary estimates.}
	We argue by contradiction. Suppose that there are no
	$\varepsilon>0$ and $\kappa>0$ for which
	\eqref{eq:quadratic_growth} holds.
	Then there exists a sequence
	$(q_n,u_n)\in W_{\mathrm{ad}}$ such that
	\[
	0<
	\|q_n-\bar q\|_{L^2(\omega)}
	+\|u_n-\bar u\|_{L^2(\Omega)}
	\le\frac1n
	\]
	and
	\begin{equation}\label{eq:sosc_contradiction}
		J(q_n,u_n)
		<
		J(\bar q,\bar u)
		+\frac{1}{2n}\left(
		\|q_n-\bar q\|_{L^2(\omega)}^2
		+\|u_n-\bar u\|_{L^2(\Omega)}^2
		\right).
	\end{equation}
	Set
	\[
	\rho_n:=
	\left(
	\|q_n-\bar q\|_{L^2(\omega)}^2
	+\|u_n-\bar u\|_{L^2(\Omega)}^2
	\right)^{1/2},
	\qquad
	h_n:=\frac{q_n-\bar q}{\rho_n},
	\qquad
	v_n:=\frac{u_n-\bar u}{\rho_n}.
	\]
	Then $\rho_n>0$, $\rho_n\to0$, and
	\begin{equation}\label{eq:sosc_normalization}
		\|h_n\|_{L^2(\omega)}^2
		+\|v_n\|_{L^2(\Omega)}^2=1.
	\end{equation}
	Consequently, after passing to a subsequence, there exist
	$h\in L^2(\omega)$ and $v\in L^2(\Omega)$ such that
	\[
	h_n\rightharpoonup h
	\quad\text{in }L^2(\omega),
	\qquad
	v_n\rightharpoonup v
	\quad\text{in }L^2(\Omega).
	\]
	
	Since $(q_n,u_n)\in W_{\mathrm{ad}}$ and $(q_n,u_n)\to(\bar q,\bar u)$ in
	$L^2(\omega)\times L^2(\Omega)$,
	Lemma~\ref{lem:admissible_state_estimate} yields a constant
	$M>0$, independent of $n$, such that
	\[
	\|u_n\|_{L^\infty(\Omega)}\le M,
	\qquad
	\|\bar u\|_{L^\infty(\Omega)}\le M.
	\]
	
For later use, set
\[
\delta u_n:=u_n-\bar u=\rho_n v_n,
\qquad
a_1(x):=\partial_t d(x,\bar u(x)),
\qquad
a_2(x):=\partial_{tt}d(x,\bar u(x)).
\]
By Assumptions~\ref{ass:d} and~\ref{ass:d_second_order},
we have $a_1,a_2\in L^\infty(\Omega)$.

We apply Taylor's formula to $d(x,\cdot)$ for almost
every fixed $x\in\Omega$.
The first-order expansion gives
\begin{equation}\label{eq:sosc_taylor_first}
	d(x,u_n(x))-d(x,\bar u(x))
	=a_1(x)\delta u_n(x)+R_n(x),
\end{equation}
where
\[
R_n(x)
=
\delta u_n(x)^2
\int_0^1(1-s)\,
\partial_{tt}d\bigl(x,\bar u(x)+s\delta u_n(x)\bigr)\,ds.
\]
Since the segment between $\bar u(x)$ and $u_n(x)$
lies in $[-M,M]$, the local uniform boundedness of
$\partial_{tt}d$ yields a constant $C_M>0$, independent
of $n$, such that
\begin{equation}\label{eq:sosc_first_order_remainder}
	|R_n(x)|\le C_M|\delta u_n(x)|^2
	\qquad\text{for almost every }x\in\Omega.
\end{equation}

\medskip
\textit{Step 2: The limit direction belongs to the critical cone.}
We show that $(h,v)\in\mathcal C(\bar q,\bar u)$.

Since $q_n\in Q_{\mathrm{ad}}$ and $\rho_n>0$, we have
\[
h_n\ge0
\quad\text{a.e. on }\{\bar q=q_a\},
\qquad
h_n\le0
\quad\text{a.e. on }\{\bar q=q_b\},
\]
with the corresponding condition omitted if a bound is infinite.
Thus, $h_n\in T_{Q_{\mathrm{ad}}}(\bar q)$.
This cone is closed and convex in $L^2(\omega)$ and hence
weakly closed. Consequently,
\[
h\in T_{Q_{\mathrm{ad}}}(\bar q).
\]

Next, we establish additional compactness for $\{v_n\}$
and verify the linearized state equation.
Subtracting the state equations for $u_n$ and $\bar u$,
testing with $\delta u_n$, and using
$\partial_t d\ge-b$, we obtain
\[
\|\nabla\delta u_n\|_{L^2(\Omega)}^2
\le
(q_n-\bar q,\delta u_n)_\omega
+b\|\delta u_n\|_{L^2(\Omega)}^2.
\]
Dividing by $\rho_n^2$ yields
\[
\|\nabla v_n\|_{L^2(\Omega)}^2
\le
\|h_n\|_{L^2(\omega)}\|v_n\|_{L^2(\Omega)}
+b\|v_n\|_{L^2(\Omega)}^2
\le \frac12+b.
\]
Thus, $\{v_n\}$ is bounded in $H^1_0(\Omega)$.
After passing to a further subsequence, the weak $L^2$ limit
$v$ from Step~1 therefore satisfies
\begin{equation}\label{eq:sosc_state_direction_convergence}
	v\in H^1_0(\Omega),
	\qquad
	v_n\rightharpoonup v
	\quad\text{in }H^1_0(\Omega),
	\qquad
	v_n\to v
	\quad\text{in }L^2(\Omega).
\end{equation}

By the first-order expansion \eqref{eq:sosc_taylor_first},
the difference of the state equations, divided by $\rho_n$,
reads
\[
(\nabla v_n,\nabla\varphi)+(a_1v_n,\varphi)
=(h_n,\varphi)_\omega
-\left(\frac{R_n}{\rho_n},\varphi\right)
\qquad\text{for all }\varphi\in H^1_0(\Omega).
\]
Since $N\in\{2,3\}$, the continuous embedding
$H^1_0(\Omega)\hookrightarrow L^4(\Omega)$ and
\eqref{eq:sosc_first_order_remainder} imply
\[
\left\|\frac{R_n}{\rho_n}\right\|_{L^2(\Omega)}
\le C_M\rho_n\|v_n\|_{L^4(\Omega)}^2
\longrightarrow0.
\]
Using $a_1\in L^\infty(\Omega)$, we may therefore pass
to the limit and obtain
\[
(\nabla v,\nabla\varphi)+(a_1v,\varphi)
=(h,\varphi)_\omega
\qquad\text{for all }\varphi\in H^1_0(\Omega).
\]
This is the linearized state equation
\eqref{eq:linearized_state}.

It remains to prove the criticality condition \eqref{eq:criticality_cond}.
The control variational inequality \eqref{eq:optimal_control} gives
\[
(\alpha\bar q+\bar z,h_n)_\omega
=\frac{1}{\rho_n}
(\alpha\bar q+\bar z,q_n-\bar q)_\omega
\ge0.
\]
To obtain the reverse inequality, we use the exact
quadratic expansion of $J$:
\[
\begin{aligned}
	J(q_n,u_n)-J(\bar q,\bar u)
	={}&\rho_n J'(\bar q,\bar u)(h_n,v_n)\\
	&+\frac{\rho_n^2}{2}
	\left(
	\alpha\|h_n\|_{L^2(\omega)}^2
	+\|v_n\|_{L^2(\Omega)}^2
	\right).
\end{aligned}
\]
Since the quadratic term is nonnegative,
\eqref{eq:sosc_contradiction} implies
\[
J'(\bar q,\bar u)(h_n,v_n)
\le
\frac{J(q_n,u_n)-J(\bar q,\bar u)}{\rho_n}
<
\frac{\rho_n}{2n}.
\]
Since $\bar q\in L^2(\omega)$ and
$\bar u-u_d\in L^2(\Omega)$, weak convergence gives
\[
\begin{aligned}
	J'(\bar q,\bar u)(h_n,v_n)
	&=\alpha(\bar q,h_n)_\omega+(\bar u-u_d,v_n)\\
	&\longrightarrow
	\alpha(\bar q,h)_\omega+(\bar u-u_d,v)\\
	&=J'(\bar q,\bar u)(h,v).
\end{aligned}
\]
Together with the preceding inequality and $\rho_n\to0$,
this implies
\[
J'(\bar q,\bar u)(h,v)\le0.
\]
On the other hand, the nonnegativity of
$(\alpha\bar q+\bar z,h_n)_\omega$ established above
passes to the limit by weak convergence.
Together with the adjoint equation
\eqref{eq:optimal_adjoint} and the linearized state
equation \eqref{eq:linearized_state}, this gives
\[
0\le(\alpha\bar q+\bar z,h)_\omega
=J'(\bar q,\bar u)(h,v)\le0.
\]
Hence the criticality condition
\eqref{eq:criticality_cond} holds, and consequently
$(h,v)\in\mathcal C(\bar q,\bar u)$.

\medskip
\textit{Step 3: Contradiction.}
By the second-order expansion from \cref{lem:cost_expansion} we have
\[
\begin{aligned}
	J(q_n,u_n)-J(\bar q,\bar u)
	={}&\rho_n(\alpha\bar q+\bar z,h_n)_\omega
	+\frac{\rho_n^2}{2}\mathcal Q(h_n,v_n) + \mathcal R(q_n,u_n)
\end{aligned}
\]
Since $(\alpha\bar q+\bar z,h_n)_\omega\ge0$,
the contradiction assumption \eqref{eq:sosc_contradiction}
implies
\begin{equation}\label{eq:sosc_quadratic_bound}
	\mathcal Q(h_n,v_n)
	\le\frac1n-\frac{2}{\rho_n^2}\mathcal R(q_n,u_n).
\end{equation}
By \eqref{eq:cost_expansion_remainder} and
$\delta u_n=\rho_n v_n$,
\[
\frac{|\mathcal R(q_n,u_n)|}{\rho_n^2}
\le
\frac{L_M}{6}\|\bar z\|_{L^\infty(\Omega)}
\rho_n\|v_n\|_{L^3(\Omega)}^3
\longrightarrow0.
\]
Here, we have used the boundedness of $\{v_n\}$ in
$H^1_0(\Omega)$ established in Step~2 and the continuous
embedding $H^1_0(\Omega)\hookrightarrow L^3(\Omega)$.
Consequently,
\begin{equation}\label{eq:sosc_quadratic_limsup}
	\limsup_{n\to\infty}\mathcal Q(h_n,v_n)\le0.
\end{equation}

Since $v_n\to v$ in $L^2(\Omega)$ and
$a_2\bar z\in L^\infty(\Omega)$, we have
\[
\int_\Omega a_2\bar z\,v_n^2\,dx
\longrightarrow
\int_\Omega a_2\bar z\,v^2\,dx.
\]
Together with the weak lower semicontinuity of the
control norm, this yields
\[
\mathcal Q(h,v)
\le\liminf_{n\to\infty}\mathcal Q(h_n,v_n)
\le0.
\]
By Step~2, $(h,v)\in\mathcal C(\bar q,\bar u)$.
The strict positivity assumption
\eqref{eq:second_order_sufficient} therefore implies
$(h,v)=(0,0)$.

In particular, $v_n\to0$ in $L^2(\Omega)$, and hence
\[
\|v_n\|_{L^2(\Omega)}^2
-\int_\Omega a_2\bar z\,v_n^2\,dx
\longrightarrow0.
\]
Using the definition of $\mathcal Q$ and
\eqref{eq:sosc_quadratic_limsup}, we conclude that
\[
0\le
\alpha\limsup_{n\to\infty}\|h_n\|_{L^2(\omega)}^2
=
\limsup_{n\to\infty}\mathcal Q(h_n,v_n)
\le0.
\]
Since $\alpha>0$, it follows that
$h_n\to0$ in $L^2(\omega)$.
Thus,
\[
\|h_n\|_{L^2(\omega)}^2+\|v_n\|_{L^2(\Omega)}^2
\longrightarrow0,
\]
contradicting \eqref{eq:sosc_normalization}.
This proves the asserted quadratic growth condition.
\end{proof}

\section{Necessary second-order optimality conditions}

In this section, we complement the sufficient optimality conditions of
Theorem~\ref{thm:second_order_sufficient} with a necessary
second-order condition.
Under Assumption~\ref{ass:inactive_set}, we show that
the quadratic form $\mathcal Q$ is nonnegative on the same critical cone
$\mathcal C(\bar q,\bar u)$.
Thus, the necessary and sufficient conditions differ
only in nonnegativity versus strict positivity of
the quadratic form on nonzero critical directions,
providing a no-gap characterization in this sense.

\begin{theorem}[Second-order necessary optimality condition]
	\label{thm:second_order_necessary}
	Let the standing assumptions on the problem data and
	Assumptions~\ref{ass:d} and~\ref{ass:d_second_order} hold.
	Let $(\bar q,\bar u)\in W_{\mathrm{ad}}$ be a local solution
	of \eqref{eq:ocp} in the sense of
	Definition~\ref{def:local_solution}, and suppose that
	Assumption~\ref{ass:inactive_set} holds.
	Let $\bar z \in H^1_0(\Omega)\cap C(\bar \Omega)$ be the unique adjoint state provided by
	Theorem~\ref{thm:optimality_system}.
	Then
	\[
	\mathcal Q(h,v)\ge0
	\qquad
	\text{for every }(h,v)\in\mathcal C(\bar q,\bar u).
	\]
\end{theorem}

We prepare the proof with three auxiliary results:
the construction of control directions compensating
for the kernel of the linearized state operator,
the construction of admissible curves for radial
critical directions, and a density result for
radial critical directions.

\begin{lemma}\label{lem:kernel_control_directions}
	Let $(\bar q,\bar u)\in W_{\mathrm{ad}}$ and suppose that
	$|I(\bar q)|>0$.
	Set $a:=\partial_t d(\cdot,\bar u)$ and define
	\begin{equation}\label{eq:kernen_K}
	K:=
	\left\{
	\psi\in H^1_0(\Omega):
	(\nabla\psi,\nabla\varphi)+(a\psi,\varphi)=0
	\quad\forall\varphi\in H^1_0(\Omega)
	\right\}.
	\end{equation}
	Then $K$ is finite-dimensional, with
	$m:=\dim K\ge0$, and
	$K\subset C(\overline\Omega)$.
	
	Moreover, there exist a measurable set $E\subset I(\bar q)$
	of positive Lebesgue measure and a number $\delta>0$ such that
	\[
	q_a+\delta\le\bar q(x)\le q_b-\delta
	\qquad\text{for almost every }x\in E,
	\]
	where inequalities involving an infinite bound are omitted.
	
	If $m>0$, let $\psi_1,\ldots,\psi_m$ be an
	$L^2(\Omega)$-orthonormal basis of $K$.
	Then there exist functions
	$e_1,\ldots,e_m\in L^\infty(\Omega)$, vanishing almost
	everywhere on $\Omega\setminus E$, such that
	\begin{equation}\label{eq:kernel_control_biorthogonality}
		(e_j,\psi_i)_\omega=\delta_{ij}
		\qquad\text{for }i,j=1,\ldots,m.
	\end{equation}
\end{lemma}

\begin{proof}
	Since $\bar u\in L^\infty(\Omega)$,
	Assumption~\ref{ass:d} implies $a\in L^\infty(\Omega)$.
	The Dirichlet realization of $-\Delta+a$ in $L^2(\Omega)$
	is self-adjoint and has compact resolvent.
	Hence its kernel $K$ is finite-dimensional.
	For every $\psi\in K$, we may rewrite its equation as
	\[
	-\Delta\psi+(a+b)\psi=b\psi.
	\]
	Since $a+b\ge0$ and $b\psi\in L^2(\Omega)$,
	Lemma~\ref{lem:monotone_equation} yields
	$\psi\in C(\overline\Omega)$.
	
	For $k\in\mathbb N$, define
	\[
	E_k:=
	\left\{
	x\in\omega:
	q_a+\frac1k\le\bar q(x)\le q_b-\frac1k
	\right\}.
	\]
	Since $I(\bar q)=\bigcup_{k\in\mathbb N}E_k$ up to
	a set of measure zero and $|I(\bar q)|>0$,
	there exists $k$ such that $|E_k|>0$.
	We choose $E:=E_k$ and $\delta:=1/k$.
	
	Suppose now that $m>0$, and introduce the Gram matrix
	\[
	G_{i\ell}:=\int_E\psi_i\psi_\ell\,dx,
	\qquad i,\ell=1,\ldots,m.
	\]
	This matrix is positive definite.
	Indeed, if $\xi\in\mathbb R^m$ satisfies
	$\xi^\top G\xi=0$, then
	\[
	\psi:=\sum_{i=1}^m\xi_i\psi_i\in K
	\]
	vanishes almost everywhere on $E$.
	By Lemma~\ref{lem:unique_continuation}, $\psi=0$ in $\Omega$.
	The linear independence of the $\psi_i$ implies $\xi=0$.
	
	We may therefore define
	\[
	e_j(x):=\chi_E(x)
	\sum_{\ell=1}^m(G^{-1})_{\ell j}\psi_\ell(x),
	\qquad x\in\Omega.
	\]
	These functions belong to $L^\infty(\Omega)$,
	vanish outside $E$, and satisfy
	\eqref{eq:kernel_control_biorthogonality}.
\end{proof}

For the next step we define the cone of radial critical directions $\mathcal C_{\mathrm{rad}}(\bar q,\bar u)$ at $(\bar q,\bar u)\in W_{\mathrm{ad}}$ as follows:
\[
\mathcal C_{\mathrm{rad}}(\bar q,\bar u)
=
\left\{
(h,v)\in\mathcal C(\bar q,\bar u):
\ \bar q+\tau h\in Q_{\mathrm{ad}}
\text{ for some }\tau>0
\right\},
\]
where $\mathcal C(\bar q,\bar u)$ is the critical cone introduced in \cref{def:critical_cone}. By convexity of $Q_{\mathrm{ad}}$, this implies
$\bar q+th\in Q_{\mathrm{ad}}$ for all $t\in[0,\tau]$.

The following lemma constructs admissible curves in $W_{\mathrm{ad}}$ with
 prescribed radial critical directions.
 
 \begin{lemma}\label{lem:admissible_curve}
 	Let the standing assumptions on the problem data and
 	Assumptions~\ref{ass:d} and~\ref{ass:d_second_order} hold.
 	Let $(\bar q,\bar u)\in W_{\mathrm{ad}}$ satisfy the
 	first-order optimality system \eqref{eq:optimality_system},
 	and suppose that $|I(\bar q)|>0$.
 	Let $K$, $m$, $E$, and $e_1,\ldots,e_m$ be as in
 	Lemma~\ref{lem:kernel_control_directions}.
 	
 	For every $(h,v)\in\mathcal C_{\mathrm{rad}}(\bar q,\bar u)$,
 	there exist $s_0>0$, $C>0$, and continuous mappings
 	\[
 	c:[0,s_0)\to\mathbb R^m,
 	\qquad
 	w:[0,s_0)\to H^1_0(\Omega)\cap C(\overline\Omega),
 	\]
 	with $c(0)=0$ and $w(0)=0$, such that
 	\[
 	q(s):=\bar q+sh+\sum_{j=1}^m c_j(s)e_j,
 	\qquad
 	u(s):=\bar u+sv+w(s)
 	\]
 	satisfy
 	\[
 	(q(s),u(s))\in W_{\mathrm{ad}}
 	\qquad\text{for all }s\in[0,s_0).
 	\]
 	Moreover,
 	\[
 	(w(s),\psi)=0
 	\qquad\text{for all }\psi\in K
 	\text{ and }s\in[0,s_0),
 	\]
 	and
 	\begin{equation}\label{eq:admissible_curve_correction}
 		|c(s)|_{\mathbb R^m}
 		+\|w(s)\|_{H^1_0(\Omega)}
 		+\|w(s)\|_{C(\overline\Omega)}
 		\le Cs^2
 		\qquad\text{for all }s\in[0,s_0).
 	\end{equation}
 	In particular,
 	\[
 	\frac{q(s)-\bar q}{s}\to h
 	\quad\text{in }L^2(\omega),
 	\qquad
 	\frac{u(s)-\bar u}{s}\to v
 	\quad\text{in }H^1_0(\Omega)\cap C(\overline\Omega)
 	\]
 	as $s\downarrow0$.
 	If $m=0$, the control correction is understood to be zero, i.e., $q(s):=\bar q+sh$.
 \end{lemma}
 \begin{proof}
 	\textit{Step 1: Preservation of the control constraints.}
 	Since $(h,v)\in\mathcal C_{\mathrm{rad}}(\bar q,\bar u)$,
 	there exists $\tau>0$ such that
 	$\bar q+\tau h\in Q_{\mathrm{ad}}$.
 	For $s\in[0,\tau/2]$ and $c\in\mathbb R^m$, define
 	\[
 	q_{s,c}:=\bar q+sh+\sum_{j=1}^m c_j e_j.
 	\]
 	By convexity $\oq+sh \in  Q_{\mathrm{ad}}$ for every $s\in[0,\tau/2]$. Moreover, the uniform distance of $\bar q$ from the
 	finite control bounds on $E$ implies
 	\[
 	\bar q+sh = \left(1-\frac{s}{\tau}\right)\bar q
 	+\frac{s}{\tau}(\bar q+\tau h) \ge \left(1-\frac{s}{\tau}\right)(q_a+\delta)
 	+\frac{s}{\tau}q_a = q_a + \left(1-\frac{s}{\tau}\right)\delta \ge q_a+\frac{\delta}{2}
 	\]
 	for almost every $x\in E$. In the same way we show the corresponding upper bound. Thus,
 	\[
 	q_a+\frac{\delta}{2}
 	\le \bar q(x)+sh(x)
 	\le q_b-\frac{\delta}{2}
 	\qquad\text{for almost every }x\in E.
 	\]
 	
 	Set
 	\[
 	C_e:=
 	\left(\sum_{j=1}^m
 	\|e_j\|_{L^\infty(\omega)}^2\right)^{1/2}.
 	\]
 	Then
 	\[
 	\left\|\sum_{j=1}^m c_j e_j\right\|_{L^\infty(\omega)}
 	\le C_e|c|_{\mathbb R^m}.
 	\]
 	Choose $r>0$ such that $C_e r\le\delta/2$.
 	Since the correction vanishes almost everywhere outside $E$,
 	we conclude that
 	\[
 	q_{s,c}\in Q_{\mathrm{ad}}
 	\qquad\text{whenever }
 	0\le s\le\tau/2
 	\text{ and }|c|_{\mathbb R^m}\le r.
 	\]
 	For $m=0$, the correction is absent and the same conclusion
 	holds without a restriction on $c$.
 	
 \medskip
 \textit{Step 2: Construction of the state and control corrections.}
 Set
 \[
 a:=\partial_t d(\cdot,\bar u)
 \]
 and let $A=-\Delta+a$ denote the Dirichlet realization
 in $L^2(\Omega)$, with domain
 \[
 D(A):=\{w\in H^1_0(\Omega):\Delta w\in L^2(\Omega)\}
 \]
 and kernel $K \subset H^1_0(\Omega)$ defined in \cref{eq:kernen_K}.
 Introduce the Banach spaces
 \[
 Y:=H^1_0(\Omega)\cap C(\overline\Omega),
 \qquad
 X:=Y\cap K^\perp,
 \]
 where orthogonality is understood in $L^2(\Omega)$.
 We equip $Y$ and its closed subspace $X$ with the norm
 \[
 \|w\|_Y
 :=\|w\|_{H^1_0(\Omega)}
 +\|w\|_{C(\overline\Omega)}.
 \]
 
 For every $y\in D(A)$, we have
 \[
 -\Delta y+(a+b)y=Ay+by.
 \]
 Since $a+b$ is bounded and nonnegative,
 Lemma~\ref{lem:monotone_equation} implies $y\in Y$ and
 \begin{equation}\label{eq:graph_norm_embedding}
 	\|y\|_Y
 	\le C\left(
 	\|Ay\|_{L^2(\Omega)}+\|y\|_{L^2(\Omega)}
 	\right).
 \end{equation}
 Thus, $D(A)$ embeds continuously into $Y$ when equipped
 with the graph norm.
 In particular, the linearized state equation
 $Av=h\chi_\omega$ implies $v\in D(A)\subset Y$.
 
 For $f\in L^2(\Omega)$, define
 \[
 Pf:=\sum_{j=1}^m(f,\psi_j)e_j,
 \]
 with $Pf=0$ if $m=0$.
 This is a bounded linear operator on $L^2(\Omega)$,
 and \eqref{eq:kernel_control_biorthogonality} yields
 \[
 (Pf-f,\psi_i)=0
 \qquad\text{for }i=1,\ldots,m.
 \]
 Hence $Pf-f\in K^\perp$.
 By the Fredholm alternative, there exists a unique
 $w\in D(A)\cap K^\perp$ satisfying
 \[
 Aw=Pf-f.
 \]
 The boundedness of the inverse of $A$ on $K^\perp$
 gives
 \[
 \|w\|_{L^2(\Omega)}
 \le C\|Pf-f\|_{L^2(\Omega)}
 \le C\|f\|_{L^2(\Omega)}.
 \]
 Together with \eqref{eq:graph_norm_embedding}, this implies
 \[
 \|w\|_Y\le C_S\|f\|_{L^2(\Omega)}.
 \]
 Consequently, the solution operator
 \[
 S:L^2(\Omega)\to X,
 \qquad A(Sf)=Pf-f,
 \]
 is linear and bounded, with a bound independent of
 $h$, $v$, and $s$.
 
We now use the Banach fixed-point theorem to construct
a correction $w(s)\in X$ and coefficients
$c(s)\in\mathbb R^m$ such that
\[
u(s):=\bar u+sv+w(s)
\]
solves the state equation
\[
\begin{aligned}
	-\Delta u(s)+d(u(s))
	&=(\bar q+sh)\chi_\omega+\sum_{j=1}^m c_j(s)e_j
	&&\text{in }\Omega,\\
	u(s)&=0
	&&\text{on }\partial\Omega.
\end{aligned}
\]
To this end, define the nonlinear remainder
\[
N(\xi):=d(\bar u+\xi)-d(\bar u)-a\xi.
\]
Since $\bar u$ satisfies the state equation and
$Av=h\chi_\omega$, the equation for $u(s)$ is equivalent to
\[
Aw(s)+N(sv+w(s))
=\sum_{j=1}^m c_j(s)e_j.
\]
Testing this equation with $\psi_i$ and using
$A\psi_i=0$ and
\eqref{eq:kernel_control_biorthogonality}, we obtain
\[
c_i(s)=(N(sv+w(s)),\psi_i),
\qquad i=1,\ldots,m.
\]
Consequently, the equation is equivalent to
\[
Aw(s)=PN(sv+w(s))-N(sv+w(s)),
\]
with the coefficients $c_i(s)$ determined by the above
identities. This motivates the fixed-point formulation below.
 Set $M:=\|\bar u\|_{L^\infty(\Omega)}+1$.
 By the local uniform boundedness of $\partial_{tt}d$,
 there exists $C_N>0$ such that
 \[
 |\partial_{tt}d(x,t)|\le C_N
 \qquad\text{for almost every }x\in\Omega
 \text{ and all }|t|\le M.
 \]
 Taylor's formula and the mean value theorem therefore give,
 for all $\xi,\eta\in L^\infty(\Omega)$ with
 $\|\xi\|_{L^\infty(\Omega)}\le1$ and
 $\|\eta\|_{L^\infty(\Omega)}\le1$,
 \begin{equation}\label{eq:curve_remainder_estimates}
 	\begin{aligned}
 		\|N(\xi)\|_{L^2(\Omega)}
 		&\le \frac{C_N}{2}
 		\|\xi\|_{L^\infty(\Omega)}\|\xi\|_{L^2(\Omega)},\\
 		\|N(\xi)-N(\eta)\|_{L^2(\Omega)}
 		&\le C_N\bigl(
 		\|\xi\|_{L^\infty(\Omega)}
 		+\|\eta\|_{L^\infty(\Omega)}
 		\bigr)\|\xi-\eta\|_{L^2(\Omega)}.
 	\end{aligned}
 \end{equation}
 The constant $C_N$ is independent of $\xi$, $\eta$,
 $h$, $v$, and $s$.
 
 For $s\ge0$, consider the fixed-point mapping
 \[
 T_s(w):=S\bigl(N(sv+w)\bigr).
 \]
 Set $V:=\|v\|_Y$.
 Combining the boundedness of $S$ with
 \eqref{eq:curve_remainder_estimates}, we obtain a constant
 $C_0>0$, independent of $h$, $v$, and $s$, such that
 \[
 \|T_s(w)\|_Y
 \le C_0(sV+\|w\|_Y)^2
 \]
 and
 \[
 \|T_s(w_1)-T_s(w_2)\|_Y
 \le C_0\bigl(
 2sV+\|w_1\|_Y+\|w_2\|_Y
 \bigr)\|w_1-w_2\|_Y,
 \]
 whenever the arguments of $N$ have $L^\infty$ norms
 at most one.
 
 Choose $R>0$ such that
 \[
 R\ge C_0(V+1)^2,
 \]
 and then choose $s_0>0$ sufficiently small that
 \[
 Rs_0\le1,
 \qquad
 s_0(V+1)\le1,
 \qquad
 2C_0s_0(V+1)\le\frac12.
 \]
 For $0\le s<s_0$, define the closed ball
 \[
 B_s:=\{w\in X:\|w\|_Y\le Rs^2\}.
 \]
 For $w\in B_s$, we have
 \[
 \|sv+w\|_{L^\infty(\Omega)}
 \le sV+Rs^2
 \le s(V+1)\le1.
 \]
 Moreover,
 \[
 \|T_s(w)\|_Y
 \le C_0s^2(V+1)^2
 \le Rs^2.
 \]
 Thus, $T_s$ maps $B_s$ into itself.
 For $w_1,w_2\in B_s$, we also have
 \[
 \|T_s(w_1)-T_s(w_2)\|_Y
 \le 2C_0s(V+1)\|w_1-w_2\|_Y
 \le\frac12\|w_1-w_2\|_Y.
 \]
 The Banach fixed-point theorem therefore yields
 $w(s)\in B_s$ satisfying
 \[
 w(s)=S\bigl(N(sv+w(s))\bigr).
 \]
 In particular,
 \[
 w(s)\perp K,
 \qquad
 \|w(s)\|_Y\le Rs^2,
 \qquad
 w(0)=0.
 \]
 
 Define
 \[
 c_j(s):=(N(sv+w(s)),\psi_j),
 \qquad j=1,\ldots,m.
 \]
 By orthonormality of the $\psi_j$,
 \[
 |c(s)|_{\mathbb R^m}
 \le \|N(sv+w(s))\|_{L^2(\Omega)}
 \le C_1s^2(V+1)^2,
 \]
 where $C_1$ is independent of $h$, $v$, and $s$.
 For $m=0$, the coefficient vector is absent.
 The fixed-point equation implies
 \[
 Aw(s)+N(sv+w(s))
 =\sum_{j=1}^m c_j(s)e_j.
 \]
 Using the state equation for $\bar u$ and
 $Av=h\chi_\omega$, we conclude that
 \[
 u(s):=\bar u+sv+w(s)
 \]
 satisfies
 \[
 -\Delta u(s)+d(u(s))
 =q(s)\chi_\omega,
 \qquad
 q(s):=\bar q+sh+\sum_{j=1}^m c_j(s)e_j|_\omega.
 \]
 After further decreasing $s_0$ so that
 \[
 s_0\le\tau/2,
 \qquad
 C_1s_0^2(V+1)^2\le r,
 \]
 Step~1 ensures $q(s)\in Q_{\mathrm{ad}}$.
 Hence $(q(s),u(s))\in W_{\mathrm{ad}}$ for all
 $s\in[0,s_0)$, and
 \[
 |c(s)|_{\mathbb R^m}+\|w(s)\|_Y
 \le \bigl(R+C_1(V+1)^2\bigr)s^2.
 \]
 This proves \eqref{eq:admissible_curve_correction}.
 
 It remains to verify continuity.
 For $s,t\in[0,s_0)$, the same remainder estimates give
 \[
 \begin{aligned}
 	\|w(s)-w(t)\|_Y
 	&\le C_0\bigl(
 	sV+\|w(s)\|_Y+tV+\|w(t)\|_Y
 	\bigr)\\
 	&\qquad\times
 	\bigl(V|s-t|+\|w(s)-w(t)\|_Y\bigr)\\
 	&\le\frac12
 	\bigl(V|s-t|+\|w(s)-w(t)\|_Y\bigr).
 \end{aligned}
 \]
 Therefore,
 \[
 \|w(s)-w(t)\|_Y\le V|s-t|,
 \]
 and $w$ is continuous.
 The definition of $c_j(s)$ and
 \eqref{eq:curve_remainder_estimates} imply continuity of $c$.
 Finally, the asserted tangent limits follow from
 \eqref{eq:admissible_curve_correction}.
 \end{proof}

 \begin{lemma}\label{lem:radial_critical_density}
 	Let Assumption~\ref{ass:d} hold, and let
 	$(\bar q,\bar u)\in W_{\mathrm{ad}}$ satisfy the
 	first-order optimality system
 	\eqref{eq:optimality_system} with an adjoint state $\bar z$.
 	Suppose that $|I(\bar q)|>0$.
 	Then, for every
 	$(h,v)\in\mathcal C(\bar q,\bar u)$, there exists a sequence
 	\[
 	(h_n,v_n)\in\mathcal C_{\mathrm{rad}}(\bar q,\bar u),
 	\qquad h_n\in L^\infty(\omega),
 	\]
 	such that
 	\[
 	\begin{aligned}
 		h_n&\to h
 		&&\text{in }L^2(\omega),\\
 		v_n&\to v
 		&&\text{in }H^1_0(\Omega)\cap C(\overline\Omega).
 	\end{aligned}
 	\]
 \end{lemma}
 
 \begin{proof}
 	Let $(h,v)\in\mathcal C(\bar q,\bar u)$ be given, and set
 	\[
 	g:=\alpha\bar q+\bar z|_\omega.
 	\]
 	Let $A=-\Delta+\partial_t d(\cdot,\bar u)$ denote the
 	Dirichlet realization in $L^2(\Omega)$.
 	We use the kernel $K$, its orthonormal basis
 	$\psi_1,\ldots,\psi_m$, and the set $E$, constant $\delta>0$,
 	and functions $e_1,\ldots,e_m$ from
 	Lemma~\ref{lem:kernel_control_directions}.
 	
 	\medskip
 	\noindent
 	\textit{Step 1: Approximation by bounded radial control directions.}
 	Following the construction in
 	\cite[proof of Theorem~7.3.5]{VexlerMeidner:2025},
 	we define
 	\[
 	\widetilde h_n(x):=
 	\begin{cases}
 		0,
 		&\text{if }q_a<\bar q(x)<q_a+\frac1n,\\
 		0,
 		&\text{if }q_b-\frac1n<\bar q(x)<q_b,\\
 		\Pi_{[-n,n]}(h(x)),
 		&\text{otherwise},
 	\end{cases}
 	\]
 	where $\Pi_{[-n,n]}$ denotes the projection onto $[-n,n]$.
 	
 	We have $\widetilde h_n\to h$ almost everywhere in $\omega$
 	and $|\widetilde h_n|\le |h|$.
 	Hence, by dominated convergence,
 	\[
 	\widetilde h_n\to h
 	\qquad\text{in }L^2(\omega).
 	\]
 	Moreover,
 	\[
 	\|\widetilde h_n\|_{L^\infty(\omega)}\le n.
 	\]
 	The variational inequality \eqref{eq:optimal_control} and
 	$h\in T_{Q_{\mathrm{ad}}}(\bar q)$, see \eqref{eq:def:T_Qad}, imply
 	$gh\ge0$ almost everywhere.
 	Since $(g,h)_\omega=0$ by criticality \eqref{eq:criticality_cond}, it follows that
 	$gh=0$ almost everywhere.
 	The definition of $\widetilde h_n$ therefore yields
 	\[
 	g\widetilde h_n=0
 	\qquad\text{almost everywhere in }\omega.
 	\]
 	
 	As in the cited proof, the truncation near the bounds,
 	the boundedness of $\widetilde h_n$, and the tangent sign
 	conditions on the active sets imply that, for every $n$,
 	there exists $\tau_n>0$ such that
 	\[
 	\bar q+s\widetilde h_n\in Q_{\mathrm{ad}}
 	\qquad\text{for all }0\le s\le\tau_n.
 	\]
 	
 	\medskip
 	\noindent
 	\textit{Step 2: Correction of the compatibility conditions.}
 	Since $(h,v)$ fulfills the linearized equation \eqref{eq:linearized_state} and $\psi_j \in K$, see \eqref{eq:kernen_K}, we have
 	\[
 	(h,\psi_j)_\omega = (\nabla v,\nabla\psi_j)
 	+(\partial_t d(\cdot,\bar u)v,\psi_j)=0,
 	\qquad j=1,\ldots,m.
 	\]
 	Define
 	\[
 	\beta_{nj}:=(\widetilde h_n-h,\psi_j)_\omega,
 	\qquad
 	r_n:=\sum_{j=1}^m\beta_{nj}e_j,
 	\qquad
 	h_n:=\widetilde h_n-r_n|_\omega.
 	\]
 	All sums are understood to be zero if $m=0$.
 	The convergence established in Step~1 implies
 	\[
 	\beta_{nj}\to0,
 	\qquad
 	\|r_n\|_{L^\infty(\Omega)}\to0.
 	\]
 	Consequently,
 	\[
 	h_n\in L^\infty(\omega),
 	\qquad
 	h_n\to h\quad\text{in }L^2(\omega).
 	\]
 	By \eqref{eq:kernel_control_biorthogonality},
 	\[
 	(h_n,\psi_i)_\omega
 	=(\widetilde h_n,\psi_i)_\omega-\beta_{ni}
 	=(h,\psi_i)_\omega=0,
 	\qquad i=1,\ldots,m.
 	\]
 	
 	The control variational inequality \eqref{eq:optimal_control} gives $g=0$
 	almost everywhere on $I(\bar q)$.
 	Since $r_n$ vanishes outside $E\subset I(\bar q)$,
 	we thus obtain
 	\[
 	gh_n=0
 	\qquad\text{almost everywhere in }\omega.
 	\]
 	
 	We next verify that $h_n$ is radial.
 	Set
 	\[
 	s_n:=
 	\min\left\{
 	\tau_n,\,
 	\frac{\delta}{1+\|h_n\|_{L^\infty(\omega)}}
 	\right\}>0.
 	\]
 	On $\omega\setminus E$, we have
 	$h_n=\widetilde h_n$, so that
 	$\bar q+s_nh_n$ satisfies the control bounds by Step~1.
 	On $E$, the distance of $\bar q$ to each finite bound
 	is at least $\delta$, whereas
 	\[
 	s_n|h_n|\le\delta
 	\qquad\text{almost everywhere}.
 	\]
 	Hence
 	\[
 	\bar q+s_nh_n\in Q_{\mathrm{ad}},
 	\]
 	which proves the radiality of $h_n$.
 	
 	\medskip
 	\noindent
 	\textit{Step 3: Approximation of the state direction.}
 	The compatibility conditions established above imply
 	\[
 	(h_n-h)\chi_\omega\in K^\perp.
 	\]
 	By the Fredholm alternative, there exists a unique
 	$w_n\in D(A)\cap K^\perp$ satisfying
 	\[
 	Aw_n=(h_n-h)\chi_\omega.
 	\]
 	The boundedness of the inverse of $A$ on $K^\perp$,
 	together with \eqref{eq:graph_norm_embedding}, gives
 	\[
 	\|w_n\|_{H^1_0(\Omega)}
 	+\|w_n\|_{C(\overline\Omega)}
 	\le C\|h_n-h\|_{L^2(\omega)}
 	\longrightarrow0.
 	\]
 	Define $v_n:=v+w_n$.
 	Since $v\in D(A)$, the graph-norm embedding also yields
 	$v\in H^1_0(\Omega)\cap C(\overline\Omega)$.
 	Thus,
 	\[
 	Av_n=h_n\chi_\omega,
 	\qquad
 	v_n\to v
 	\quad\text{in }H^1_0(\Omega)\cap C(\overline\Omega).
 	\]
 	Finally, using the adjoint equation, we obtain
 	\[
 	\begin{aligned}
 		J'(\bar q,\bar u)(h_n,v_n)
 		&=\alpha(\bar q,h_n)_\omega
 		+(\bar u-u_d,v_n)\\
 		&=(g,h_n)_\omega=0.
 	\end{aligned}
 	\]
 	Together with the radiality of $h_n$, this proves
 	\[
 	(h_n,v_n)\in\mathcal C_{\mathrm{rad}}(\bar q,\bar u)
 	\]
 	and completes the proof.
 \end{proof}
 
 \begin{proof}[Proof of Theorem~\ref{thm:second_order_necessary}]
 	\textit{Step 1: Nonnegativity on radial critical directions.}
 	Let $(h,v)\in\mathcal C_{\mathrm{rad}}(\bar q,\bar u)$.
 	By Lemma~\ref{lem:admissible_curve}, there exist admissible
 	curves of the form
 	\[
 	q(s)=\bar q+sh+k(s),
 	\qquad
 	u(s)=\bar u+sv+w(s),
 	\]
 	where
 	\[
 	k(s):=\sum_{j=1}^m c_j(s)e_j|_\omega
 	\]
 	and
 	\[
 	\|k(s)\|_{L^2(\omega)}
 	+\|w(s)\|_{H^1_0(\Omega)}
 	+\|w(s)\|_{C(\overline\Omega)}
 	\le Cs^2.
 	\]
 	Set
 	\[
 	\delta q(s):=q(s)-\bar q,
 	\qquad
 	\delta u(s):=u(s)-\bar u.
 	\]
 	Then
 	\[
 	\begin{aligned}
 		\frac{\delta q(s)}{s}&\to h
 		&&\text{in }L^2(\omega),\\
 		\frac{\delta u(s)}{s}&\to v
 		&&\text{in }H^1_0(\Omega)\cap C(\overline\Omega)
 	\end{aligned}
 	\]
 	as $s\downarrow0$.
 	In particular,
 	\begin{equation}\label{eq:bound_du_s}
 	\|\delta u(s)\|_{L^\infty(\Omega)}
 	+\|\delta u(s)\|_{H^1_0(\Omega)}
 	\le Cs.
 	\end{equation}
 	
 	Set $g:=\alpha\bar q+\bar z|_\omega$.
 	The control variational inequality
 	\eqref{eq:optimal_control} implies that $g=0$
 	almost everywhere on $I(\bar q)$.
 	Since the functions $e_j$ vanish outside
 	$E\subset I(\bar q)$, we have $(g,k(s))_\omega=0$.
 	Together with the criticality condition
 	\eqref{eq:criticality_cond}, this gives
 	\begin{equation}\label{eq:nosc_control_cancellation}
 		(g,\delta q(s))_\omega
 		=s(g,h)_\omega+(g,k(s))_\omega=0.
 	\end{equation}

 	Since $(\bar q,\bar u)$ is a local solution, we may
 	decrease $s_0$ such that
 	\begin{equation}\label{eq:nosc_local_comparison}
 		0\le J(q(s),u(s))-J(\bar q,\bar u)
 		\qquad\text{for all }s\in[0,s_0).
 	\end{equation}
 	
By \eqref{eq:bound_du_s}, there exists $M>0$,
independent of $s\in[0,s_0)$, such that
 	\[
 	\|u(s)\|_{L^\infty(\Omega)} \le M \quad\text{and} \quad \|\bar u\|_{L^\infty(\Omega)} \le M
 	\]
 	Using the quadratic expansion from \cref{lem:cost_expansion} we obtain
 	\[
 	\begin{aligned}
 	0\le J(q(s),u(s))-J(\bar q,\bar u) &= (\alpha\bar q+\bar z,\delta q(s))_\omega
 	+\frac12\mathcal Q(\delta q(s),\delta u(s))
 	+\mathcal R(q(s),u(s))\\
 	&=\frac12\mathcal Q(\delta q(s),\delta u(s))
 	+\mathcal R(q(s),u(s)),
 	\end{aligned}
 	\]
 	where we also used \eqref{eq:nosc_control_cancellation}.
 	
 	For $0<s<s_0$ we obtain by \cref{lem:cost_expansion}
 	\[
 	\begin{aligned}
 	\mathcal Q\left(\frac{1}{s}\delta q(s),\frac{1}{s}\delta u(s)\right) &\ge -\frac{2}{s^2} \mathcal R(q(s),u(s))\\
 	& \ge -\frac{L_M}{3s^2} \|\bar z\|_{L^\infty(\Omega)} \|\delta u(s)\|_{L^3(\Omega)}^3\\
 	& \ge-\frac{C}{s^2} \|\bar z\|_{L^\infty(\Omega)} \|\delta u(s)\|_{H^1_0(\Omega)}^3 \ge -Cs,
 	\end{aligned}
 	\]
 	where we also used the embedding $H^1_0(\Omega) \hookrightarrow L^3(\Omega)$ and \eqref{eq:bound_du_s}. 	
 	
	By the convergence of the difference quotients and
	the continuity of $\mathcal Q$ on
	$L^2(\omega)\times L^2(\Omega)$,
	passing to the limit as $s\downarrow0$ yields
 	\[
 	\mathcal Q(h,v)\ge0
 	\qquad
 	\text{for every }
 	(h,v)\in\mathcal C_{\mathrm{rad}}(\bar q,\bar u).
 	\]
 	
 	\medskip
 	\noindent
 	\textit{Step 2: Extension to all critical directions.}
 	Let $(h,v)\in\mathcal C(\bar q,\bar u)$.
 	By Lemma~\ref{lem:radial_critical_density}, there exists
 	a sequence
 	\[
 	(h_n,v_n)\in\mathcal C_{\mathrm{rad}}(\bar q,\bar u)
 	\]
 	such that
 	\[
 	h_n\to h\quad\text{in }L^2(\omega),
 	\qquad
 	v_n\to v\quad\text{in }L^2(\Omega).
 	\]
 	Step~1 and the continuity of $\mathcal Q$ imply
 	\[
 	\mathcal Q(h,v)
 	=\lim_{n\to\infty}\mathcal Q(h_n,v_n)\ge0.
 	\]
 	This completes the proof.
 \end{proof}
 
 \section*{Acknowledgments}
 The author used ChatGPT (OpenAI) to assist with literature
 searches, review of mathematical
 arguments, and proofreading. The author assumes responsibility for all content.
 
\bibliographystyle{siamplain}
\bibliography{references}
\end{document}